\documentclass[msom,nonblindrev]{informs3} 
\OneAndAHalfSpacedXII 

\usepackage{placeins}
\usepackage{pifont}
\usepackage[T1]{fontenc}
\usepackage[utf8]{inputenc}
\usepackage{authblk}
\usepackage{graphicx}
\usepackage{multicol}
\usepackage[sectionbib,round]{natbib}
\title {Adjustable Robust Machine Scheduling}
\date{}
\author{}
\usepackage{url}
\usepackage[backref=page,unicode,naturalnames,colorlinks=false]{hyperref}
\usepackage[dvipsnames]{xcolor}
\usepackage{empheq} 
\usepackage{ulem}
\usepackage{float}
\usepackage{tikz}
\usepackage{accents}
\newcommand{\ubar}[1]{\underaccent{\bar}{#1}}
\usepackage[shortlabels]{enumitem}
\usetikzlibrary{matrix,decorations.pathreplacing, calc, positioning}
\usepackage[english]{babel}
\usepackage{caption}
\usepackage{subcaption}
\usepackage{comment}
\mathtoolsset{showonlyrefs=true} %shows equation number if it is referenced in the text
\usepackage{booktabs}
\usepackage{algorithm}
\usepackage{algorithmic}
\usepackage{multirow}
\usepackage{siunitx}
\usepackage{fancyhdr} 
\newcommand{\ie}{{\it i.e.}}%,\xspace}
\newcommand{\eg}{{\it e.g.}}%,\xspace}
\TheoremsNumberedThrough     % Preferred (Theorem 1, Lemma 1, Theorem 2)
\newtheorem{prop}{Proposition}
\EquationsNumberedThrough    % Default: (1), (2), ...
\usepackage{glossaries}
\usepackage{array}
\usepackage{makecell}
\newacronym[plural=RCPSPs,firstplural=Resource-Constrained Project Scheduling Problems (RCPSPs)]{RCPSP}{RCPSP}{Resource-Constrained Project Scheduling Problem}
\newacronym{MILO}{MILO}{Mixed Integer Linear Optimization}
\newacronym{RO}{RO}{Robust Optimization}
\newacronym{ARO}{ARO}{Adaptive Robust Optimization}
\newacronym{2SSA}{2SSA}{Two-stage Static Allocation}
\newacronym{SA}{SA}{Static Allocation}
\newacronym{SL}{SL}{Static List}
\newacronym{AR}{AR}{Adaptive Robust}
\newacronym{PMS}{PMS}{Parallel Machine Scheduling}
\newacronym{PH}{PH}{Perfect Hindsight}
\newacronym{DP}{DP}{Dynamic Programming}
\newacronym{RH}{RH}{Rolling Horizon}
\newacronym{CDF}{CDF}{Cumulative Distribution Function}

\newcommand{\revise}[1]{\textcolor{black}{#1}}
 \DeclareMathOperator*{\pinf}{inf\vphantom{p}}          
\newcommand{\cmark}{\ding{51}}
\newcommand{\xmark}{\ding{55}}
\newcommand{\hcmark}{\ding{51}{\raisebox{0.1em}{\kern-0.6em\footnotesize\ding{55}}}}

\begin{document}
%%%%%%%%%%%%%%%%

%\setcounter{page}{1}
%\YEAR{0000}% \eg,  2005
%\DOI{10.1287/xxxx.0000.0000}%

\RUNAUTHOR{}

%\RUNTITLE{Robust Machine Scheduling}

% \TITLE{Bundling Information Goods of Decreasing Value}
\TITLE{An Adaptive Robust Optimization Model for Parallel Machine Scheduling}

\ARTICLEAUTHORS{%
\AUTHOR{Izack Cohen}
\AFF{Faculty of Engineering, Bar-Ilan University, Ramat Gan, Israel, \EMAIL{ izack.cohen@biu.ac.il}}
\AUTHOR{Krzysztof Postek}
\AFF{Faculty of Electrical Enginnering, Mathematics and Computer Science, Delft University of Technology, Delft, The Netherlands, \EMAIL{k.s.postek@tudelft.nl}}
\AUTHOR{Shimrit Shtern}
\AFF{Faculty of Industrial Engineering and Management, Technion - Israel Institute of Technology, Haifa, Israel, \EMAIL{shimrits@technion.ac.il}}
} % end of the block

% Short abstract
\ABSTRACT{Real-life parallel machine scheduling problems can be characterized by: (i) limited information about the exact task duration at scheduling time, and (ii) an opportunity to reschedule the remaining tasks each time a task processing is completed and a machine becomes idle. Robust optimization is the natural methodology to cope with the first characteristic of duration uncertainty, yet the existing literature on robust scheduling does not explicitly consider the second characteristic – the possibility to adjust decisions as more information about the tasks’ duration becomes available, despite that re-optimizing the schedule every time new information emerges is standard practice. In this paper, we develop a scheduling approach that takes into account, at the beginning of the planning horizon, the possibility that scheduling decisions can be adjusted. We demonstrate that the suggested approach can lead to better here-and-now decisions and better makespan guarantees. To that end, we develop the first mixed integer linear programming model for adjustable robust scheduling, and a scalable two-stage approximation heuristic, where we minimize the worst-case makespan. Using this model, we show via a numerical study that adjustable scheduling leads to solutions with better and more stable makespan realizations compared to static approaches.}


\maketitle
\section{Introduction} \label{sec:intro}
\gls{PMS} problems are widely researched owing to their theoretical importance and multiple applications in manufacturing, cloud computing, and project management, among others. Real-life \gls{PMS} settings involve uncertainty about task duration, which may be characterized by the randomness of each task duration and, possibly, a dependence between task durations.

An ideal scheduling approach should accommodate uncertainty to ensure realistic guarantees on the objective function value and permit adjustments of later-stage scheduling decisions based on different observed task lengths (\eg,  different duration realizations of the task scheduled first may result in different allocation decisions of the next tasks).

In the existing literature, commonly, scheduling problems are analyzed using various types of \textit{static schedules}. \gls{SA} is one of the most frequently used static policies, where the tasks are allocated in a certain order to pre-specified machines. In such policies, both the decision about the order of tasks on each machine and their allocation to machines are static, and do not change as more information is revealed. Another example are \gls{SL} policies, where the scheduling order of tasks to machines is pre-specified and tasks are allocated, by this order, to the first idle machine. These policies can be viewed as having a static order of task allocations but adaptive in the choice of the machine to which each task is allocated. For both types of policies, research has addressed the issue of finding the optimal policy that minimizes the expectation/worst-case value of the objective function over uncertain task durations. For both policy types, however, it holds that the optimal \gls{SA} and \gls{SL} policies may change once the duration of some completed tasks is known. We are not aware of works that consider this issue explicitly.

\begin{comment}
Before attempting to bridge this gap in the treatment of static policies, we should decide on the way the decision maker (hereafter, \textit{scheduler}) accounts for the uncertainty when making scheduling decisions. Specifically, we need to decide what is the performance measure to optimize and whether to optimize for its expected value or worst-case value.
\end{comment}

We focus on a makespan minimization objective function that is used for load balancing in \gls{PMS} and many other scheduling applications. When deciding whether to use the expected value or worst-case value, several factors should be considered. Optimizing over an expectation requires specifying the full probability distribution of task duration, information that is often not readily available or is costly to acquire. Moreover, the makespan of a single realization may significantly differ from the expected value; thus, if the exact scheduling problem is not repeated multiple times, optimizing over the expected value may not be translated into good performance in practice. In contrast, much less information is needed when specifying a set that includes all the reasonable duration realizations, and a worst-case optimization approach provides a guarantee on the performance of any realization in such a set. Therefore, we consider a scheduler who minimizes the worst-possible makespan of a set of tasks over some uncertainty set, which captures all reasonable scenarios within the support of the distribution. This is in line with the paradigm of \gls{RO}, where the best solution is sought under the assumption that the problem's parameters are initially unknown and that, given the decisions, nature picks their worst-possible values from an \textit{uncertainty set} consisting of outcomes that include the true realization with a high probability. A representative real-world \gls{PMS} example is presented by \cite{xu2013robust} who describe a new product development division that needs to manufacture prototypes for multiple newly developed mobile phones. These prototypes include sets of new parts that are being manufactured for the first time. In the absence of historical data regarding processing duration, stochastic models are irrelevant and \gls{RO} becomes a leading alternative for hedging against schedule delays.

In this mindset, we consider the classical version of \gls{PMS}, where $m$ identical machines process $n \geqslant m$ tasks that are available at the start of the scheduling horizon. We construct an exact \gls{MILO} formulation for minimizing the worst-case makespan, which includes all possible later-stage (re-)scheduling decisions and gives the best-possible \gls{AR} policy. Since this formulation scales exponentially in the problem size, we also propose a scalable adaptive heuristic -- the \gls{2SSA} --  where only one re-optimization moment is considered. Next, we compare the adaptive formulation optimal scheduling decisions and optimal worst-case makespan to those of the optimal \gls{SA} and \gls{SL}.

In contrast to the majority of previous works, which compare naive implementations of the \gls{SA} and \gls{SL} policies without re-optimization (\ie, re-scheduling) as more information is revealed, we consider the more realistic rolling horizon implementation of these policies. Under this implementation, whenever one of the machines becomes idle, the scheduler can alter the initial order of tasks by re-solving an optimization problem with the extra information included. \revise{A summary of the considered policies, in decreasing order of the solution quality in our experiments, is given in the following table:
\begin{center}
\begin{tabular}{lcc}
 & Adaptive & Scalable \\ \hline
\gls{AR}: Adaptive & \cmark & \xmark \\ \hline
\gls{2SSA}: Two-stage static & \hcmark & \cmark \\ \hline
\gls{SL}: Static list & \hcmark & \xmark \\ \hline
\gls{SA}: Static allocation & \xmark & \cmark \\ \hline
\end{tabular}
\end{center}}

\vspace{0.5cm}

The main contributions of our research are as follows:
\begin{enumerate}
\item \revise{We characterize settings in which using adaptive policies may be important. Specifically, we identify settings in which the static and adaptive policies result in the same  makespan, and provide performance bounds for the \gls{SA} policy as well as any rolling horizon policy with respect to a \gls{PH} policy, which determines the best allocation when the durations are known exactly. The bounds are computed for the popular {\it budgeted uncertainty set} \citep{bertsimas2004price}. Using these bounds we determine that when the budget is moderate relative to the number of tasks $n$, and when $n$ is not extremely large, planning for adaptivity may be crucial for good promised and actual performance.}

\item We provide a closed-form \gls{MILO} formulation of the \gls{PMS} problem with re-scheduling and its scalable heuristic counterpart. To the best of our knowledge, this is the first such formulation. \revise{The scalable heuristic counterpart, named \gls{2SSA}, accounts for one stage of adaptivity. We demonstrate, via experiments, the scalability of \gls{2SSA} and the benefit it derives from partially accounting for adaptivity in large problems.}
%\item We investigate the adaptive policies and identify a setting in which the static and adaptive policies result inpartially the same \shimrit{or in a similar} makespan.
\item We demonstrate, through stylised examples and numerical experiments, that both the adaptive robust policy and its scalable heuristic counterpart can significantly outperform their static alternatives, even when the latter are implemented via a rolling-horizon approach.
\item In terms of managerial insights, the main conclusion of our paper is that, when possible, future re-scheduling of tasks (adaptations) should be taken into account at the planning stage as a way to obtain substantially better makespan guarantees as well as to shorter actual makespans compared to non-adaptive approaches.
\end{enumerate}	

The remainder of the paper is structured as follows. In Section~\ref{sec:literature.review}, we review the relevant scheduling and robust optimization literature. Section~\ref{sec:notation} introduces the notation and definitions used in our formulations and analyses. In Section~\ref{sec:StaticVsAdjustable}, we discuss structural properties of static and  adaptive policies. \revise{In Section \ref{sec:Dynamic-programming}, we introduce the \gls{DP} formulations of the adaptive robust scheduling from the scheduler's and adversary's points of view. In Section \ref{sec:MIO} we develop, via the adversary view, the \gls{MILO} formulation of the problem, as well as its scalable \gls{2SSA} counterpart.} Section \ref{sec:numerical.study} demonstrates, through a numerical study, the benefit of using adaptive policies. Section~\ref{sec:insights} presents managerial insights gained from our investigation of adaptive policies and Section~\ref{sec:future.research} concludes and suggests future research directions.

\section{Literature review} \label{sec:literature.review}
\begin{comment}
This paper focuses on adaptive robust makespan minimization for the \gls{PMS} problem. Task durations are assumed to be within a so-called bounded uncertainty set, and allocation decisions can be adjusted as new information emerges. Makespan minimization, and hence balancing the load across machines, is an important real-world objective. The classical, deterministic version of the \gls{PMS} with identical machines, which is one of the most studied \gls{PMS} problems \citep{ranjbar2012two}, is NP-Hard \citep{pinedo2000scheduling}.
\end{comment}
We focus on adaptive robust makespan minimization for the classical \gls{PMS} problem with identical machines. The deterministic version of the problem, which is one of the most studied \gls{PMS} problems \citep{ranjbar2012two}, is NP-Hard \citep{pinedo2000scheduling}. We review the two main approaches for dealing with uncertain durations in the context of \gls{PMS}: the stochastic optimization approach and the robust optimization approach. The former approach is attractive if probability distributions of task durations are known and the scheduler desires a policy that performs well on average (see, Section \ref{lit:stochastic}). If, in contrast, the scheduler wants assurance that the policy will perform well for any realization of the durations within a predefined set, or does not have an accurate estimate of the underlying distribution, then a robust (min-max) approach that minimizes the worst-case performance is the best option (see, Sections \ref{lit:static}-\ref{lit:adjustable}).

Although this paper does not address the stochastic setting, due to the connection between the stochastic and robust settings, our literature review will first cover the more studied stochastic models, before transitioning to discuss the \gls{RO} setting and its adaptive variant. % in this context.

\subsection{Stochastic \texorpdfstring{\gls{PMS}}{PMS} models} \label{lit:stochastic}
This type of scheduling model optimizes an expected value objective, such as the expected makespan to process $n$ tasks on $m$ machines. The probability distributions of task durations are assumed to be known or can be inferred based on historical data.

There are only a few known optimal policies for specific probability distributions. For example, the longest expected processing time (LEPT) priority rule -- by which tasks are processed according to a non-increasing order of their expected duration -- minimizes the expected makespan for exponentially distributed and independent task durations \citep{cai2014optimal}.

Many studies looked at the stochastic \gls{PMS} under various conditions and assumptions. A partial list of representative publications includes: \cite{mohring1999approximation} who developed non-anticipative scheduling policies via linear programming relaxations to minimize the expected weighted flow time (that is, the sum of expected task completion times). They analyzed the performance of the weighted shortest expected processing time priority rule, which is simple to apply, and found that it is asymptotically optimal; \cite{ranjbar2012two} developed efficient branch-and-bound procedures to maximize the probability that a set of tasks with normally-distributed processing times completes before its due date. Their experiments included up to 20 tasks and five identical machines; \cite{weber1982scheduling} allowed preemptions and developed priority rules based on highest/lowest hazard rates; others solved \gls{PMS} problems using heuristics such as genetic algorithms and simulated annealing \citep[\eg, ][]{balin2011parallel,yeh2014parallel}.

Commonly, when historical data are not accessible or costly to acquire or when the tasks are unique, probability distributions of task durations are not available and schedulers may have to rely on estimations of upper and lower bounds of task durations \citep{balouka2019robust}. For non-repetitive tasks, decision makers tend to exhibit a risk-averse behavior that hedges against worst-case scenarios \citep{daniels1995robust, lin2011robust}. In such cases, the use of \gls{RO}, which we review next, is natural.  

\subsection{Static robust PMS}\label{lit:static}
To the best of our knowledge, all previous research applying \gls{RO} to the \gls{PMS} problem involved static policies that do not consider, when making scheduling decisions, the possibility that decisions could or should be changed later on. We suggest the option to adjust task/machine allocations as new information is uncovered. Below are several studies that use static robust solution approaches.

The closest work to ours is by \cite{xu2013robust} who investigated the robust \gls{PMS} with identical machines and a makespan minimization objective under processing times specified via an interval-type uncertainty. The authors formulated the problem as a \gls{MILO}, and solved it via exact solution approaches based on iterative relaxation algorithms and several heuristics. A computational experiment with up to five machines and 15 tasks demonstrated that the heuristics' average deviation from the optimal value is smaller than 8\%. The current research departs from \cite{xu2013robust} by developing an adjustable \gls{RO} model, which may use convex uncertainty sets including (but not limited to) the conservative interval-type uncertainty set, showing the equivalence between the static and adaptive solution in this case. \cite{bougeret2019robust} focused on developing approximation algorithms for robust \gls{PMS} with identical machines and a budgeted uncertainty set. 

Other studies adopted a min-max regret objective function that minimizes the maximal deviation of a given solution from the optimum across all scenarios, since such objective is considered as less conservative than the traditional robust objective \citep{aissi2009min}. \cite{conde2014mip} formulated a \gls{MILO} problem to minimize the maximal regret of the flow time  for a \gls{PMS} environment with unrelated machines in which the processing durations are specified via an interval-type uncertainty set. They solved the \gls{MILO} for problems with up to 40 tasks and 10 machines using a bound on the computation time. Importantly, even a simpler version of this problem with identical machines is NP-Hard \citep{de2016impact}. \cite{xu2014hedge} who investigated \gls{PMS} with unrelated machines showed that a solution with the nominal (midpoint) task processing durations is a two-approximation for the min-max regret problem. They suggested an interesting modeling idea by which the original problem is transformed into a robust single machine problem, which yielded an \gls{MILO} problem with fewer variables and constraints. The authors stressed the importance of researching \gls{PMS} problems with other objective functions -- an idea we adopt by using the makespan minimization objective. 

Other papers addressed robust \gls{PMS} problems with a variety of objectives and constraints such as cost minimization where task processing can be outsourced \citep{wang2020approximation} and maximization of the probability that all tasks complete by their due dates using a distributionally robust approach \citep{liu2019service}.

\subsection{Adjustable robust optimization}\label{lit:adjustable}
Most \gls{RO} research in scheduling implements a `static' approach, as detailed in the previous section. Yet, in problems spanning multiple time periods, most schedulers would adjust their decisions as the actual duration of each task emerges. A common approach for emulating this is the rolling-horizon approach, where the problem is re-solved at selected decision points, taking into account the new information. The downside of this approach is that if a static \gls{RO} approach is used to solve each of the respective optimization problems, the here-and-now decisions are still optimized as if the later-stage decisions were not to going to change regardless of the uncertainty realization. 

A branch of \gls{RO} research dealing with this issue is \gls{ARO} \citep{ben2004adjustable}. Its main idea is to optimize the here-and-now decisions, taking into account all scenarios in which the problem may unfold and the corresponding optimal decisions for each such scenario. Optimal \gls{ARO} solutions are favorable with respect to their static counterparts since they result in better here-and-now decisions that take into account adaptations in later time stages. At the same time, solving problems via \gls{ARO} is computationally demanding because of the need to account for the contingent decisions in a large number of scenarios. A particularly difficult case is problems in which later-stage decisions are discrete. 

Since solving \gls{ARO} exactly is often computationally demanding, one may restrict the space of considered policies to obtain a tractable approximation. The most common restriction is using affine decision rules, as introduced by \cite{ben2004adjustable}. There, later-stage decisions are affine functions of the unknown parameters, and the coefficients of these functions are optimized as decision variables. For some problems, affine decision rules are shown to be optimal \citep{bertsimas2010optimality}, but this is not true in general. Modeling of continuous variables as affine decision rules may yield computationally efficient solutions as \cite{cohen2007stochastic} demonstrate for the time-cost tradeoff project scheduling problem. 

The four main approaches to approximating \gls{ARO} in the case of integer decisions, which we deal with in scheduling, are: (i) the $K$-adaptability approach of \cite{hanasusanto2015k}, (ii) the iterative partitioning approach of \cite{postek2016multistage} and \cite{bertsimas2016multistage}, (iii) cutting-plane-like approaches as in \cite{zeng2013solving}, and (iv) decision rule approximations \citep[see,][and references therein]{georghiou2019decision}. All these approaches assume that one knows the time moments at which the values of initially unknown parameters are revealed. This makes them applicable to problems such as unit commitment \citep{bertsimas2012adaptive}, inventory control \citep{ben2004adjustable}, or flood protection planning \citep{postek2019adjustable}. 
In the \gls{PMS} problem that we focus on, the above mentioned assumption is no longer valid. The time moments at which new information (completion of tasks) becomes available depend on (i) the uncertain durations of the currently running tasks and (ii) previous scheduling decisions. Moreover, all the decisions are discrete schedule-or-not binaries. For this reason, the \gls{PMS} problem structure is uncharted territory for \gls{ARO}, and exactly where our research makes a notable contribution.

\section{Notation and definitions} \label{sec:notation}
We consider minimizing the processing makespan of tasks $i=1,2,\dots,n$ on $m$ identical machines. Tasks and machines are available at time $t=0$, there are no precedence relations between tasks, and tasks cannot be preempted while processing. The task durations vector $d=(d_1,\ldots,d_n)$ is uncertain and lies within an \textit{uncertainty set} $U$. The scheduler aims to minimize the worst-case makespan over this predefined uncertainty set, \ie, assuming that for all scheduling decisions possible, the adversary (nature) will pick the worst-possible duration vector $d \in U$.

\begin{figure}[H]
		\caption{A schematic representation of an example timeline for a system with two machines}
		\label{fig:StatesIllustration}  
		\begin{center}	
			\includegraphics[width=0.8\textwidth]{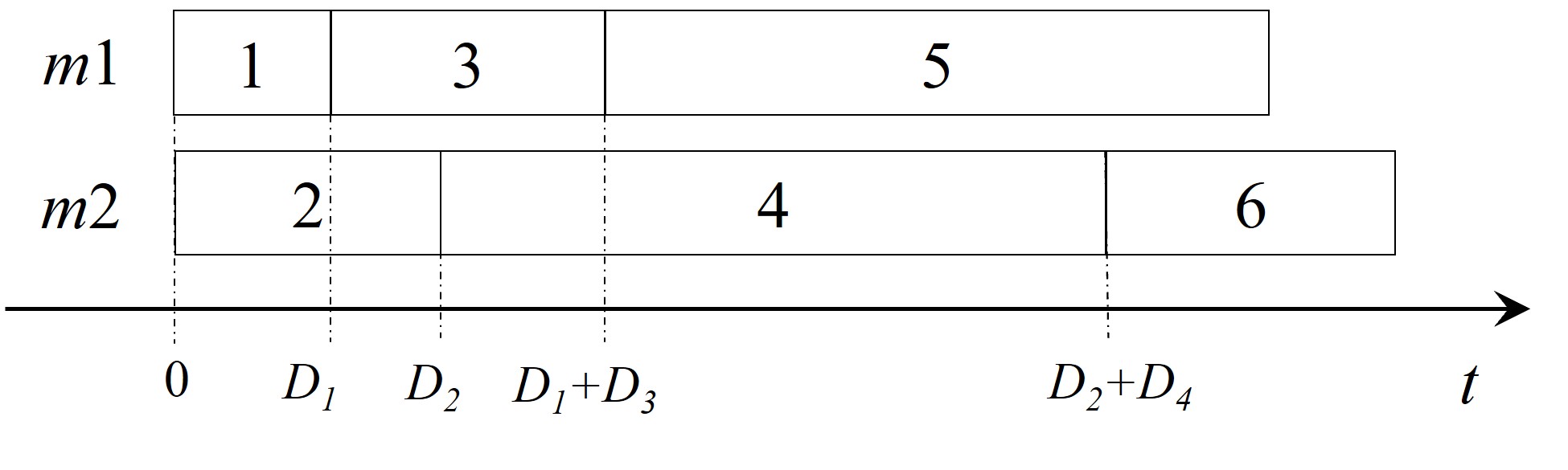}
		\end{center}
	\end{figure}

To model the problem, we need to describe the possible states of the system. Using the notation summarized in Table~\ref{table:Notation}, in what follows, we develop our modeling framework. To explain the notation for the system state consider a system with $m = 2$ machines, each of which is either processing a task or idle. When a machine completes a task and becomes idle, an unscheduled task, if one exists, has to be immediately allocated to the machine. The system state at such a time moment $t$ is described by $(S, F,D,i, \bar{D_i})$, where $S$ is an ordered list of started tasks (in the order of starting with arbitrary tie-breaks), $F$ is an ordered set of completed tasks (in the order of finishing with arbitrary tie-breaks), $D$ is an ordered list of realized durations corresponding to the completed tasks in $F$, $i$ is a task that is currently being processed and its duration $d_i$ is known to be $d_i \geqslant \bar{D_i}$ time units, where $\bar{D_i}$ is its processing time at the instance the state is observed. There are also, however, states in which both machines are busy, and when this occurs, no decision is made. Thus, to keep the state space limited, we only include states in which at least one machine is idle.  %\shimrit{[suggestion: moving this to after a policy is defined]} 
Figure~\ref{fig:StatesIllustration} and Table~\ref{table:StatesIllustration} demonstrate the states of a system with two machines for a certain realization and schedule.

As the scheduling progresses in time, the up-to-now (total) duration of the running (completed) tasks becomes known, respectively. This new information might reduce the uncertainty about the possible duration of the remaining tasks, eliminating certain parts of the uncertainty set. Therefore, we introduce the notion of state-dependent uncertainty. This notion is formally expressed in the definition of the uncertainty set induced by a system state as
$$
U_{S,F,D,i,\bar{D}_{i}}=\{d\in U: d_k = D_{k},\ \forall k\in F, \ d_{i} \geq \bar{D}_{i}\}.
$$
We call a state  $(S,F,D,i,\bar{D}_{i})$ feasible if $U_{S,F,D,i,\bar{D}_{i}} \neq \emptyset$. We define the set of feasible states as $\mathcal{S}$.

Note that the notation discussed so far for the case of $m=2$ can be extended to $m>2$ by replacing $i$ and $\bar{D_i}$ with $I$ and $\bar D$ -- the sets of in-process tasks and their processing times until $t$, respectively.

\FloatBarrier
	\begin{table}
		\caption{Primary notation}
		\label{table:Notation}  
		\begin{center}	
			\small
			{\renewcommand{\arraystretch}{1.15} %donne la distance entre les lignes%
				{\setlength{\tabcolsep}{0.2cm}
					\begin{tabular}{| c | p{10cm} |}
						\hline
						\textbf{Indices, parameters and variables}&\textbf{Description}\\ \hline
						
						$i, j, k$ & Denote a task or a machine, by context\\ \hline
						$m$ & Number of identical machines\\ \hline
						$n$ & Number of tasks\\ \hline
						$t, \tau$ & Denote time from start of processing \\ \hline
						$d_i, d$ & Duration of task $i$ and the vector of task durations, respectively\\ \hline
						$D_i, \bar{D_i}$ & The realized duration of a completed task $i \in F$, and the amount of time task $i\in I$ has been processing at the observation time $t$, respectively\\ \hline
						$s_i$ & The start time of task $i$ \\ \hline
						\textbf{Sets and lists}& \\ \hline
						$\mathbb{R}^n _{\geqslant 0}$ & An $n$-dimensional set of non-negative, real numbers \\ \hline
						$[n]$ & A shorthand for the set $\{1,\ldots,n\}$\\
						\hline
						%$V$ & Set of tasks\\ \hline
						$S$ & An ordered list of started tasks \\ \hline
						$F, \bar{D}$ & An ordered list of completed tasks and their realized durations, respectively\\ \hline
						$I, D$ & Sets of the in-process tasks at time $t$ and their respective duration till $t$ \\ \hline
						$U$ & Uncertainty set for task durations \\ \hline
						
			\end{tabular}}}
		\end{center}	
	\end{table}

We now define the notion of a scheduling \textit{policy}. We define a policy $P$ as a mapping from a state to the choice of the next task to be scheduled. That is, $P:\mathcal{S}\rightarrow [n]$, where $P(S,F,D,I,\bar{D})\in[n]\setminus S$. In this work, we compare three types of policies for minimizing the worst-case makespan. 
\begin{itemize}
\item \gls{SA} policies, where each task is assigned to a machine according to a predefined order. Thus, given an uncertainty set $U$, a \gls{SA} policy amounts to a partition of the tasks to machines $J=(J_1,\ldots,J_m)$ such that $\cup_{j\in[m]} J_j=\revise{[n]}$ and $J_j \cap J_k=\emptyset$ for all $j\neq k$. Once the decision about the partition $J$ is made, it does not change. Assuming that the tasks in each $J_j$ are given in a certain order (without loss of generality, we assume it is lexicographic), the policy can be explicitly defined as
$$P^{\text{SA},J}(S,F,D,I,\bar{D})= \argmin\{k\in[n]\setminus S: \exists j\in[m],\; k\in J_j,\; i\notin J_j \;\forall i\in I\}.$$
Thus, at each decision state, the next task to be scheduled is the first task that has not yet started and has been allocated to the machine that just became idle.
Note that this policy is independent of the information gained in each stage about the durations of the tasks that have already begun processing, meaning $D$ and $\bar{D}$, and thus, it does not adapt to this information.
\item \gls{SL} policies, where a task is processed on the first idle machine according to its location in an ordered list. Thus, given an uncertainty set $U$, the \gls{SL} policy amounts to a permutation of the tasks:
$\pi=(\pi_1,\pi_2\ldots, \pi_n)$, where for each $i_{th}$ place in the list, $\pi_i\in[n]$ is a specific task, and all tasks must be allocated once, \ie,  $\pi_i\neq \pi_j$ for any $i\neq j\in[n]$.  Once the decision on the permutation is made, it does not change. Given permutation $\pi$, the policy can be explicitly defined as
$$P^{\text{SL},\pi}(S,F,D,I,\bar{D})= \pi_{i}, \quad i=\argmin\{k\in [n]: \ \pi_k\notin S \}.$$
Thus, at each decision state, the next in order on the list of not-started tasks is scheduled on the idle machine.
Note that similarly to the \gls{SA} policy, the \gls{SL} policy also does not adapt to the state-dependent information.
\item The \gls{AR} policy, denoted by $P^{\text{AR}}$ is the most flexible and considers all possible scenarios for task realizations, choosing the next task according to the actual scenario that was realized.  In particular, the scheduling decision for two distinct states with the same completed and processing tasks may be different. Specifically,  $(D,\bar{D})\neq (D^\dag,\bar{D}^\dag)$ allows for
$P^{\text{AR}}(S,F,D,I,\bar{D})\neq P^{\text{AR}}(S,F,D^\dag,I,\bar{D}^\dag)$. 
\end{itemize}

All three policy types can be applied in a rolling-horizon fashion by solving a new optimization problem once a machine becomes idle, where the uncertainty set is the one induced by the revealed state of the system. Nevertheless, only the \gls{AR} policy takes such implementation into account at $t = 0$. Hence, we expect the optimal \gls{AR} policy to yield favorable makespan guarantees due to its improved here-and-now scheduling decisions compared to the optimal \gls{SL} and \gls{SA} policies that do not consider adaptations in later-stage decisions.
	\begin{table}[H]
		\caption{System states for the timeline presented in Figure (\ref{fig:StatesIllustration})}
		\label{table:StatesIllustration}  
		\begin{center}	
			\small
			{\renewcommand{\arraystretch}{1.15} %donne la distance entre les lignes%
				{\setlength{\tabcolsep}{0.2cm}
					\begin{tabular}{| c | p{9.7cm} |p{2.3cm} |}
						\hline
						\textbf{Time $t$}&\textbf{State notation}&\textbf{Decision}\\ \hline				
						0 & ($[], [], [],0,0$) & Start 1 and 2 \\ \hline
						$D_1$ & ($[1,2],[1],D_1,2,D_1$) & Start 3 \\ \hline
						$D_2$ & ($[1,2,3],[1,2],[D_1, D_2],3,D_2-D_1$) & Start 4 \\ \hline
						$D_1+D_3$ & ($[1,2,3,4],[1,2,3],[D_1,D_2,D_3],4,D_1+D_3-D_2$) & Start 5 \\ \hline
						$D_2+D_4$ & ($[1,2,3,4,5],[1,2,3,4],[D_1, D_2,D_3,D_4],5,D_2+D_4-(D_1+D_3)$) & Start 6 \\ \hline
			\end{tabular}}}
		\end{center}	
	\end{table}
Finally, to formulate our problem of minimizing the worst-case makespan, we introduce the function $T(S,F,D,i,\bar{D_i})$ denoting the minimal worst-case duration to process the unfinished tasks (\ie, $i$ and the tasks that have not started), assuming an optimal \gls{AR} policy is applied. We refer to this function as the \textit{remaining makespan}. Our objective of minimizing the worst-case makespan is equivalent to the function $T([],[],[],0,0)$.

In Section~\ref{sec:Dynamic-programming}, we use the above notation in developing \gls{DP} formulations for finding the optimal \gls{AR} policy.

\section{Comparing the static and adjustable robust policies} \label{sec:StaticVsAdjustable} 
In this section, we characterize and highlight potential differences between the three scheduling policies defined in the previous section. In Section \ref{sec:Box uncertainty set}, we demonstrate that if the uncertain task durations are independent of each other, the best worst-case makespan is identical for all policy types. In Section \ref{sec:Diff-RO-ARO}, we use toy examples with task duration dependence for showing that an optimal \gls{AR} policy may yield different first-stage decisions from those of its static counterparts, \gls{SA} and \gls{SL}. Based on this result, in the remainder of the paper, we focus on situations where task durations are connected to each other via the shape of the uncertainty set. Accordingly, in Section \ref{sec:bound} we provide performance bounds for the case of two machines under the commonly used budgeted uncertainty set, for which the possible tasks' durations are connected through the budget. Specifically, we bound the performance of the optimal \gls{SA} policy and of a rolling-horizon implementation of any scheduling policy.
These bounds quantify the possible gains from utilizing an adaptive policy as a function of the problem's parameters such as the number of tasks, their durations and the uncertainty set's parameters.

\subsection{Equivalence of the scheduling policies under box uncertainty sets} \label{sec:Box uncertainty set}
In Section \ref{sec:notation}, we defined a scheduling policy $P$ as a mapping from the system's states into scheduling decisions. The schedule is created by deciding which tasks to schedule and observing which task completed first and at what time. Thus, for a given policy $P$ and vectors of durations $d$, the \textit{schedule} can be represented as an ordered partition of the tasks to machines $J=(J_1,\ldots, J_m)$ where its makespan is given by $\max_{j\in[m]}\sum_{i\in J_j} d_i$. Let $\mathcal{M}(P,d)$ denote the function that maps a policy and a set of durations into such a partition. 

The optimal \gls{SA} policy is equivalent to the optimal partitioning of tasks to machines $J^*=(J^*_1,\ldots,J^*_m)$. Thus, for a given duration vector $d$, we have $\mathcal{M}(P^{\text{SA},J^*},d)\equiv J^*$, which is independent of $d$. Therefore, denoting all sets of partitions of $n$ tasks to $m$ machines as $\mathcal{J}_{n.m}$, an optimal \gls{SA} policy is obtained by solving the following optimization problem:
\begin{equation}\label{eq:static}
\min_{J\in\mathcal{J}_{n,m}}\max_{d\in \mathcal{U}}\max_{j\in[m]}\sum_{i\in J_j} d_i.
\end{equation}
In contrast, the optimal \gls{SL} policy is equivalent to an optimal permutation $\pi^*$ of the tasks. Consequently, the allocation of tasks in $\pi^*$ to machines and the respective makespan depend on the task duration vector $d$. For illustration, consider two machines, three tasks and the permutation $(1,2,3)$. The schedule will start with $1$ and $2$ on the two machines; then, if $1$ completes first, $3$ will be scheduled on its machine and otherwise on the other machine. The worst-case makespan for a \textit{given} \gls{SL} policy associated with permutation $\pi$ is, therefore:
\begin{equation}\label{eq:wc_list}
\max_{d\in\mathcal{U},J=\mathcal{M}(P^{\text{SL},\pi},d)} \max_{j\in[m]} \sum_{i\in J_j} d_i.
\end{equation}
We denote $\Pi_n$ as the set of all permutations of $[n]$. Thus, the optimal \gls{SL} policy is obtained by the solution of  \eqref{eq:wc_list} over $\Pi_n$: 
\begin{equation}\label{eq:min_wc_list}
	\min_{\pi\in \Pi_n}\max_{d\in\mathcal{U}, J=\mathcal{M}(P^{\text{SL},\pi},d)} \max_{j\in[m]} \sum_{i\in J_j} d_i.
\end{equation}
Finally, denoting the space of all possible policies as $\mathcal{P}$, the optimal \gls{AR} policy is given by the solution of the following optimization problem.
\begin{equation}\label{eq:min_wc_AR}
	\min_{P\in\mathcal{P}}\max_{d\in\mathcal{U}, J=\mathcal{M}(P,d)} \max_{j\in[m]} \sum_{i\in J_j} d_i.
\end{equation}
The box uncertainty set, by which the duration of each task is bounded within an interval and independent of the other tasks, is a natural candidate for modeling task durations in scheduling settings. This situation underlines the importance of our next result -- that the optimal makespan achieved by all three robust policies is equal when using a box uncertainty set. In other words, the here-and-now decisions produced by all three policies are optimal (for brevity, we placed all the proofs in Appendix \ref{app:proofs}). 
\begin{prop} \label{prop.box.the.same}
Let $\mathcal{U}$ be a box uncertainty set given by 
$$
\mathcal{U}=\left\{d: d_i\in [\ubar{d}_i,\bar{d}_i], i \in [n]\right\}.
$$
Then, the optimal \gls{SA} policy, the optimal \gls{SL} policy, and the optimal \gls{AR} policy produce the same worst-case makespan.
\end{prop}

\subsection{Different first-stage decisions in \texorpdfstring{\gls{RO}}{RO} and \texorpdfstring{\gls{AR}}{AR}}  \label{sec:Diff-RO-ARO} 
Following the results of the previous section, we focus on uncertainty sets that are not box shaped. This section demonstrates the superiority of the adjustable \gls{AR} policy with respect to the alternative \gls{SA} and \gls{SL} policies, even if they are implemented in a rolling-horizon fashion; that is, only the here-and-now decisions are implemented and the policies are re-optimized every time a task finishes. Rolling horizon implementation of scheduling policies is a common management practice for adapting to new information. Specifically, we show via examples that the optimal first-stage decisions of an optimal \gls{AR} policy may be different than those of optimal \gls{SA} and \gls{SL} policies, which implies better performance by the former. %For brevity, we will refer to some optimal \gls{SA}, \gls{SL}, and \gls{AR} policies as the \gls{SA}, \gls{SL}, and \gls{AR} policies, respectively.

We begin with an example with three tasks in which the optimal here-and-now decisions are the same for the \gls{SL} and \gls{AR} policies (thus these two policies are equivalent), but are different for the \gls{SA} policy. Then, we present an example with four tasks, where the \gls{AR} policy has different here-and-now optimal decisions from the \gls{SL} policy. The takeaway is that an \gls{AR} policy may achieve better objective function values compared to \gls{SA} and \gls{SL} policies (by more than 6\% in our toy examples). 
\subsubsection{A three-task example.}
For the case of three tasks, we can obtain closed-form short expressions for the worst-case makespan for each policy. Consider first the \gls{SA} policy in which, without loss of generality, tasks 1 and 2 are processed first, and task 3 is processed after task 1. Thus, $J_1=(1,3), J_2=(2)$. The worst-case duration is then
$$
\max\limits_{d \in U} \max \{ d_1  + d_3, d_2 \},
$$
which is the maximum of these two terms:
$$
\max\limits_{d \in U}d_1  + d_3 , \qquad \max\limits_{d \in U} d_2.
$$
For \gls{AR}, consider, without loss of generality, a schedule in which tasks 1 and 2 are processed first, and task 3 starts processing as soon as the first of the two tasks is finished. The worst-case duration under \gls{AR} is
$$
\max\limits_{d \in U} \max \{ \min\{d_1, d_2 \} + d_3, \max\{d_1, d_2 \} \},
$$
which is the maximum of these three terms:
$$
\max\limits_{d \in U} \min\{d_1 + d_3, d_2  + d_3\}, \qquad \max\limits_{d \in U}  d_1, \qquad \max\limits_{d\in U} d_2
.$$
To choose the optimal allocation for each policy, we can evaluate the above expressions for every permutation of the three tasks. With respect to \gls{SL}, we note that for any setting in which $n=m+1$ ($n=3$ and $m=2$ for our example), the optimal \gls{SL} policy and the optimal \gls{AR} policy are equivalent.

To make things concrete, we consider a specific setting in which task durations are:
$$
d_1 = 0.0580 + 0.95z_1 , \ d_2 = 0.1945 + 0.75 z_2, \ d_3 = 0.5866 + 0.48 z_3 ,
$$
where the uncertain parameters are $(z_1, z_2, z_3) \in Z = \{ [0, 1]^3, \sum_{i=1}^3 z_i \leq 2.5 \}$ (\ie, a budgeted uncertainty set). Let us compare the optimal \gls{AR} (here equivalent to \gls{SL}) with the static policy of \gls{SA}. The optimal solutions, respectively, are:
\begin{itemize}
\item \gls{AR}: Start with tasks 1 and 2, and start processing task 3 whenever the first of these tasks has completed. This gives a worst-case duration of 1.83.
\item \gls{SA}: There are two equivalent solutions. (i) Start tasks 1 and 3, and process task 2 after 1 or, (ii) start tasks 2 and 3, and process task 1 after 2. Both options give a worst-case duration of 1.95. Any other option in which the first tasks are 2 and 3 leads to a longer duration.
\end{itemize}
Thus, even in such a simple setting, the adjustable policy leads to a makespan lower by about 6\% compared to \gls{SA}. Moreover, since the first-stage decisions are different (\ie, `allocate $\{1,2\}$' by \gls{AR} compared to `allocate $\{1,3\}$ or $\{2,3\}$' by \gls{SA}), a rolling-horizon policy applied for the static policy will still be inferior to \gls{AR}. In fact, if we start with tasks $\{1,3\}$ and schedule task 2 after the first of them is completed, we will obtain a worst-case makespan of 1.95 (identical to the one in which we did not adapt the decision). If we start with $\{2,3\}$ and schedule task 2 after the first of them is completed, we obtain a worst-case makespan of 1.88, which is still almost 3\% more than that of \gls{AR}. Indeed, our numerical study presented in Section~\ref{sec:numerical.study} indicates that \gls{AR} may be better than \gls{SA} by up to 30\%.

\subsubsection{A four-task example.}	
 We now illustrate the difference between the three policies, \gls{SA}, \gls{SL} and \gls{AR}, for a setting with $n=4$ and $m=2$. We use an uncertainty set that includes the five scenarios outlined in Table \ref{table:ScenariosForDiscretUSet}.
\begin{table}[H]
	\caption{Specification of the uncertainty set}
	\label{table:ScenariosForDiscretUSet}  
	\begin{center}	
		\small
		{\renewcommand{\arraystretch}{1.15} %donne la distance entre les lignes%
			{\setlength{\tabcolsep}{0.2cm}
				\begin{tabular}{ c  p{2cm} p{2cm} p{2cm} p{2cm} p{2cm} }
					\hline
					{Task}&\multicolumn{5}{ c }{Task durations for:}\\ 		
					&{Sc. 1}&{Sc. 2}&{Sc. 3}&{Sc. 4}&{Sc. 5}\\ \hline				
					1& 3 & 4.5 & 4.75 & 2.5 & 0.25\\ \hline					
					2 & 2 & 2 & 2 & 3.5 & 5 \\ \hline
					3 & 3 & 3.5 & 3 & 3& 3.5 \\ \hline
					4 & 5.5 & 4 & 4 & 4 & 4\\ \hline
		\end{tabular}}}
	\end{center}	
\end{table}
Without providing the (tedious) calculations, we summarize the optimal choices of each of the three policies.
\begin{itemize}
	\item \gls{AR}: The optimal \gls{AR} policy schedules tasks 1 and 4 first. Then, for scenarios 1,2, and 3, task 3 is scheduled to start when the first of 1 and 4 completes, and task 2 is scheduled last. For scenarios 4 and 5, task 2 is scheduled to start immediately when the first of 1 and 4 completes and task 3 is last. The corresponding optimal worst-case makespan value is 7.5.
	\item \gls{SL}: There are four optimal list policies by which the minimal makespan is 8, longer by $6.7\%$ compared to the optimal \gls{AR}. The policies are: $(1,2,4,3), (1,4,2,3), (2,3,4,1),(2,4,1,3)$. Take for example, $(2,3,4,1)$. It achieves makespans of $7.5,8,7.75,7,7.5$ for scenarios $1$ to $5$, respectively. Thus, the optimal worst-case makespan is 8.
	\item \gls{SA}: The static allocation policy partitions the tasks to machines upfront. The optimal partition is $J_1=(1,2),J_2=(3,4)$. It achieves makespans of $8.5,7.5,7,7,7.5$ for scenarios  $1$ to $5$, respectively. Thus, the robust makespan is 8.5. All other partitions lead to longer makespans. The static allocation policy leads to a makespan that is longer by 13.3\% compared to \gls{AR}.
\end{itemize}

This example underscores the importance of \gls{AR} in scheduling for gaining better promised makespan guarantees and better actual makespan. Regarding the former aspect, the three policies \gls{AR}, \gls{SL}, \gls{SA} give different `promised worst-case' makespans, which is important from a managerial standpoint of providing guarantees (\eg, when submitting a contract proposal). Regarding the latter aspect, the \gls{AR} policy will perform better than both \gls{SL} and \gls{SA} even if the latter two were re-optimized every time a task finishes (rolling horizon). The difference lies in the possibility of making a sub-optimal here-and-now decision in the case of \gls{SL} and \gls{SA}. Our numerical study, in Section \ref{sec:numerical.study}, clearly shows that the performance of a policy deteriorates as the number of its different first-stage decisions with respect to \gls{AR} increases. 

\begin{comment}
{For both \gls{SL} and \gls{SA}, initial scheduling of tasks 2 and 3 is no different from scheduling of 1 and 4 (both correspond to policies with promised optimal makespans of 8 and 8.5, respectively). Nevertheless, scheduling tasks 2 and 3 first is not an optimal choice and even if re-optimization was applied at every decision state, the realized makespan will be at least 8.}
\end{comment}
\subsection{\revise{Performance bounds under the budgeted uncertainty set}} \label{sec:bound}

\revise{In this section, we identify the cases in which taking into account decision adaptivity in the planning stage may lead to a significant makespan improvement compared to static policies and rolling horizon implementations of policies. % compared to non-adaptive alternatives, such as the common \gls{SA}, and cases in which rolling horizon implementations both adaptive and non-adaptive policies may perform well. 
We focus on the two-machine case, and establish upper performance bounds for the \gls{SA} and for policies which do not allow leaving a machine idle. These latter family of policies include \gls{SL} and a \gls{RH} implementation of \gls{SA}, and are denoted by RH in the formulas below. We compare these policies to the perfect hindsight policy, which we denote as \gls{PH} -- the optimal policy when all the uncertain task durations are assumed to be known in advance. In other words, we wish to bound, from above, the ratios $\frac{\text{SA}}{\text{PH}}$ and $\frac{\text{RH}}{\text{PH}}$. The first ratio provides a bound on the maximum suboptimality of the `promised' makespan and the second one bounds the maximum suboptimality of the actual rolling-horizon `execution' of a policy. These two bounds will provide the potential benefit that can be obtained by using \gls{AR}, and help us to identify regions in which accounting for later adaptivity would be the most useful.}

\revise{To provide bounds, we need to choose an uncertainty setting. We choose the well known and widely used budgeted uncertainty set, first suggested by \cite{bertsimas2004price}. The structure of the uncertainty set is captured by the following assumption. 
\begin{assumption}\label{ass:shape_uncertainty}
The nominal length of any task $i$ is bounded, \ie,  $d^0_i\in[\ubar{a},\bar{a}]$ for some $0<\ubar{a}\leq \bar{a}<\infty$. Moreover, there exists $\alpha>0$ such that $\bar{d}_i=\alpha d^0_i$ for all $i\in[n]$, and 
$$U=\left\{d\in R^n: d_i=d_i^0+\bar{d}_i u_i, u_i\in[0,1], \sum_{i=1}^n u_i\leq \Gamma\right\}.$$
\end{assumption}
This uncertainty set implies that in the considered realizations the nominal duration vector $d^0$ can be augmented by a maximal perturbation $\bar{d}$ proportional to $d^0$ (can be thought of as a percent of the nominal). However, not all tasks will achieve this maximal perturbation, since the sum of the ratio of the perturbations used for each task cannot exceed a given budget $\Gamma$.}
%For both ratios, \shimrit{we lower bound the worst case makespan of \gls{PH}, by assuming that we can perfectly balance the two machines, thus arriving to the bound:
%$$PH\geq \sup\limits_{d \in U} \frac{\sum\limits_{i = 1}^n d_i}{2}.
%$$}

\subsubsection{\revise{Bound on the promised durations of SA versus PH.}} \label{sec:bound.promised.durations}
\revise{We start by bounding the ratio between \gls{SA} and \gls{PH}. For the case of two machines, an allocation can be represented as a binary vector $x \in \{0, 1 \}^n$  such that
$$
x_i = \left\{ \begin{array}{ll} 1 & \text{if the $i$-th task is on machine 1}, \\
0 & \text{otherwise.} \end{array} \right.
$$
With this notation, the processing time on machine 1 is $x^\top d$ and the processing time on machine 2 is $(e - x)^\top d$ where $e$ is a vector of ones. Thus, the total makespan is $\max \{ x^\top d, (e - x)^\top d  \}$.
Mathematically, the makespans of \gls{SA} and \gls{PH} are 
\begin{align*}
SA=\min\limits_{x \in \{0, 1\} } \sup\limits_{d \in U} \max \{ x^\top d, (e - x)^\top d \},
\end{align*}
and 
\begin{align*}
PH=\sup\limits_{d \in U} \min\limits_{x \in \{0, 1\} } \max \{ x^\top d, (e- x)^\top d \}.
\end{align*}
To bound the $\frac{\text{SA}}{\text{PH}}$ ratio, we will upper bound \gls{SA} and lower bound \gls{PH}. We choose a common allocation to work with, corresponding to the optimal allocation of \gls{PH} for $d=d^0$. This allocation, which might be sub-optimal for \gls{SA}, provides an upper bound on \gls{SA} that can be found by computing the worst case $d\in U$. To bound the value of this worst case, we note that the uncertainty budget will be allocated first to tasks on a single machine before allocating it to the tasks on the other machine.}

\revise{In order to find a lower bound on \gls{PH}, we use a specific choice of $d\in U$, in which the budget is allocated equally between the tasks, such that each task is perturbed by $\Gamma/n$ of its maximal perturbation. Finally, we use the fact that the optimal allocation $x$ for this realization is equal to the nominal allocation.}

\revise{Below, we summarize the bound obtained by using this approach.}
%by selecting a specific static allocation -- the one that corresponds to the nominal task durations -- and computing \shimrit{the worst case makespan for that allocation, corresponding to a specific worst case perturbation of $d$. This worst case perturbation will reflect the fact that in the worst case the budget will first be allocated to tasks on a single machine, and if $\Gamma$ is big enough, then the remaining budget will be allocated to tasks on the other machine.} To that end, we need to upper bound the maximum number of tasks on a single machine. 
%We also lower bound the duration of \gls{PH} given by the following 
%We do so by restricting the adversary to allocate an equal share of the budget, $\Gamma / n$, to each of the tasks' perturbation.

%Before proving the performance bound in Theorem \ref{theo:SA}, we need to define, via Assumption \ref{ass:shape_uncertainty}, finite upper and lower bounds on the nominal task durations and a relationship between a nominal duration and its maximal deviation. Notice, that this assumption can be applied to multiple settings.

\revise{\begin{theorem} \label{theo:SA}
Let the number of machines be $m=2$ and let Assumption~\ref{ass:shape_uncertainty} hold. Let $\tilde{x}$ be the partition to machines which minimizes the makespan for the deterministic problem where $d=d^0$. Specifically, we denote by $\mathcal{I}$ and $\bar{\mathcal{I}}=[n]\setminus \mathcal{I}$ the set of tasks that $\tilde{x}$ allocates to the first and second machine, respectively. Furthermore, we denote the ordering of $d^0_i$ for tasks in $\mathcal{I}$ and $[n]\setminus\mathcal{I}$ by
$d^{0,\mathcal{I}}_{(1)}\geq \ldots \geq d^{0,\mathcal{I}}_{(|\mathcal{I}|)}$ and 
$d^{0,\bar{\mathcal{I}}}_{(1)}\geq \ldots \geq d^{0,\bar{\mathcal{I}}}_{(n-|\mathcal{I}|)}$,  respectively.
Then, the ratio between the worst case makespan resulting from the SA policy and the PH policy is bounded from above by
\begin{equation}
   \frac{\text{SA}}{\text{PH}}\leq \frac{n}{n+\Gamma\alpha} \left ( 1 +\alpha\max_{\mathcal{S}\in \{\mathcal{I},\bar{\mathcal{I}}\}}\frac{\sum_{k=1}^{\min\{|\mathcal{S}|,\lfloor{\Gamma}\rfloor\}} d^{0,\mathcal{S}}_{(k)}+\Delta^\mathcal{S}(\Gamma)d^{0,\mathcal{S}}_{\left(\min\{|\mathcal{S}|,\lfloor\Gamma\rfloor\}+1\right)}}{\sum_{i\in\mathcal{S}} d_i^0} \right ),
\end{equation}
where for any $\mathcal{S}\subseteq[n]$ we have $\Delta^{\mathcal{S}}(\Gamma)=\min\{\Gamma,|\mathcal{S}|\}-\min\{\lfloor{\Gamma}\rfloor,|\mathcal{S}|\}$.
\end{theorem}
Indeed, the above theorem shows that for $\Gamma=0$ and $\Gamma=n$ the ratio is bounded from above by $1$, \ie, the value of adaptivity decreases as the budget goes to $0$ or $n$. Thus, for any value of $d^0$ there exists some {$0<\Gamma<n$} for which this ratio is maximized, and so is the potential value of adaptivity. We also note that as $n$ and $\Gamma$ stays proportional to $n$, the bound does not necessarily go to $1$, indicating that employing the \gls{SA} without re-optimizing may be extremely suboptimal even for large $n$.} %\color{red}DO WE WANT TO INCLUDE THIS? IT IS TRUE BUT NOT SURE IT CONTRIBUTES. We further observe that when $n\rightarrow \infty$, then for a fixed $\Gamma$, $\frac{SA}{PH}\rightarrow 1$.\color{black} 

\subsubsection{\revise{Bound on the rolling horizon implementation of any scheduling policy}}\label{sec:bound.rh}
\revise{We now turn to bound the ratio of any \gls{RH} policy and \gls{PH}. In this case, we lower bound the worst case makespan of \gls{PH}, by assuming that we can perfectly balance the two machines, thus arriving to the bound:
$$PH\geq \sup\limits_{d \in U} \frac{\sum\limits_{i = 1}^n d_i}{2}.$$}

\revise{In order to bound the \gls{RH} implementation of any scheduling policy, we make the following observation. The difference between the individual makespan of the two machines must be less than or equal to the maximum duration of the last task to finish (note that in an extreme case, this can be an extremely long task, which occupies the entire makespan of one of the machines). The identity of this last task depends on both the policy and the adversary's decision. Regardless, we can bound the makespan of any \gls{RH} implementation by bounding from above the time at which the last task starts by and adding to it the duration of the last task.
For example, if the last task to start is $i$, the time at which it will start will not be longer than $\sum_{j = 1, \ j \neq i}^n d_j/{2}$. Therefore, we can bound the \gls{RH} makespan by
$$
RH\leq \max_{1 \leq i \leq n} \sup\limits_{d \in U} \frac{\sum\limits_{j = 1, \ j \neq i}^n d_j}{2} + d_i.
$$}

\revise{Next, we introduce Theorem \ref{theo:RH}, which uses the above bounds to bound the desired ratio.
\begin{theorem}\label{theo:RH}
Let the number of machines be $m=2$ and let Assumption~\ref{ass:shape_uncertainty} hold. Then the ratio between the worst case makespan resulting from the RH policy and the PH policy is bounded from above by
\begin{equation}
   \frac{\text{RH}}{\text{PH}}\leq 1+\frac{\bar{a}\min\{(1+\alpha),(1+\alpha\Gamma)\}}{\ubar{a}(n+\alpha\Gamma)}.
\end{equation}
\end{theorem}
The above theorem implies that for any fixed $0<\ubar{a}\leq \bar{a}<\infty,\alpha>0$ and any $\Gamma>0$, as $n\rightarrow\infty$ any \gls{RH} implementation gets closer to that of \gls{PH}. Thus, for very large values of $n$ there will be almost no benefit to computing \gls{AR}, since one can get good performance by using any arbitrary \gls{RH} policy.}

\section{Dynamic programming formulations} \label{sec:Dynamic-programming} 
In the previous section, we highlighted potential differences between the three types of policies \gls{SA}, \gls{SL}, and \gls{AR}. {We now begin to develop} a general method to optimize the \gls{AR} policy. A natural first step is to formulate a scheduler-perspective \gls{DP} model of the problem -- this is the focus of this section. Solving this \gls{DP} will be computationally intractable, but its formulation lays the foundation of the adversary-perspective \gls{DP} of Section~\ref{sec:adversary.DP} that we solve via an \gls{MILO} formulation in Section~\ref{sec:MIO}. To keep the exposition clear, we consider the two-machine setting here, with the general $m$-machine cases considered in Appendices~\ref{app:DP_more_than_two} and \ref{app:DP_more_than_two_advers}.
\subsection{Scheduler's dynamic programming formulation}
At decision points, which occur at (i) the beginning of the planning horizon and (ii) when a task completes, the scheduler allocates to an idle machine a task that minimizes the worst-case makespan from that time point and on (hereafter, remaining makespan). We denote by $T(\cdot)$ the function that outputs the worst-case remaining makespan at a state described by $S,F,D,i,\bar{D}_{i}$ in which task $i$ is still being processed on a machine and the second machine is idle (\eg,  it just completed processing a task). The recursive formulation for $T$ is then given by
%{ \small
\begin{align}\label{eq:T_no_Waitng_eps}
T(S,F,D,i,\bar{D}_{i})=&\min_{k\notin S}%F\cup\{i\}}
\max \Biggl\{ \Biggr. \max_{\substack{d_k:d\in U_{[S,k], F,D,i,\bar{D}_{i}},\\d_k \leq d_i-\bar{D}_{i}}} d_k+ 
T([S,k], [F,k],[D,d_k],i,\bar{D}_{i}+d_k),\\
&\qquad\quad\max_{\substack{d_i:d\in U_{[S,k], F,D,i,\bar{D}_{i}},\\d_k>d_i-\bar{D}_{i}}} d_i-\bar{D}_{i}+T([S,k], [F,i],[D,d_i],k,d_i-\bar{D}_{i})\Biggl.\Biggr\}.\nonumber
\end{align}%}\\
The objective (outer $\min$) is to allocate a task $k$ that minimizes the worst-case remaining makespan. For each $k$, there are two possible scenarios from which an adversary chooses the worst of the first $\max$. The first term within the parentheses relates to the possibility that under the worst-case scenario, $k$ completes its processing before the completion of the processed task $i$. In such a case, the next decision point is due when $k$ completes and $i$ is still running. The second term considers the possibility that under the worst-case scenario, $i$ completes before $k$, in which case the next decision point is when task $i$ finishes processing.\footnote{In this and later formulations, we seemingly ignore the case where two or more tasks finish processing at exactly the same time. This case can be dealt with explicitly; however, it significantly complicates the presentation of the problems and is, therefore, omitted for clarity in all following formulations.}

For the boundary case, given by $|S|=n$, $|F|<n$, when all tasks have already been scheduled, we have that
\begin{align*}T(S,F,D,i,\bar{D}_{i})=\max_{d_i:d\in U_{S,F,D,i,\bar{D}_{i}}} d_i-\bar{D}_i,\end{align*}
\ie, the remaining makespan will only be the time until the current task $i$ is completed.

The takeaway from this section is that in order to solve the scheduler-perspective \glspl{DP}, we may have to optimize over an infinite dimensional space of policies, \ie, all functions from states to scheduling decisions. For that reason, determining the optimal policy is a difficult task, calling for an approach different from the usual \gls{DP} solution techniques, which we introduce in the next section.

\subsection{The adversary dynamic programming formulation}\label{sec:adversary.DP}
Our approach to determining the optimal schedule takes the perspective of an adversary who seeks the scenario with the worst-possible task durations, taking into account that the scheduler will make the best-possible scheduling decisions at each decision point. In other words, the adversary tries to make the makespan as long as possible given that the scheduler implements an optimal scheduling policy.

%In light of the above reasoning, Section~\ref{sec:adversary.DP} develops the \gls{DP} from the perspective of the adversary and Section~\ref{sec:MIO}, formulates this \gls{DP} using a \gls{MILO}. In this \gls{MILO} formulation, we incorporate all possible choices made by the adversary regarding which of the currently processing tasks to end first, as well as the scheduler decision regarding which of the tasks to schedule next, using auxiliary binary variables. The main decision variables in the \gls{MILO} are the adversary's choice of durations $d$. The optimal scheduling policy will follow implicitly from the optimal solution of this problem.

The adversary's decision points are the states in which a certain task has just been scheduled so that either (i) all machines are busy, or (ii) there are no tasks left to schedule. An example of such a state is:
$([S,k],F,D,i ,\bar{D}_{i})$,
where task $k\notin S\equiv F\cup \{i\}$ has just been scheduled. To present the recursion, we define the function $T'(\cdot)$ for all the adversary's decision points as the worst-case remaining makespan. Thus, the recursive formula for $T'$ is given by
\begin{align}\label{eq:adversary_DP_2_machines}
\begin{split}
& T'([S,k],F,D,i ,\bar{D}_{i})= \\
& \max\Biggl\{\Biggr. \max_{\substack{d_k:d\in U_{[S,k], F,D,i,\bar{D}_{i}},\\d_k \leq d_i-\bar{D}_{i}}} d_k+ \min_{l \notin S\cup\{k\}} T'([S,k,l], [F,k],[D,d_k],i,\bar{D}_{i}+d_k),\\
&\qquad\max_{\substack{d_i:d\in U_{[S,k], F,D,i,\bar{D}_{i}},\\d_k>d_i-\bar{D}_{i}}} d_i - \bar{D}_{i} + \min_{l \notin S\cup\{k\}} T'([S,k,l], [F,i],[D,d_i],k,d_i-\bar{D}_{i})\Biggl. \Biggr\},
\end{split}
\end{align}
as long as $|S| \leq  n-2$, \ie, not all tasks have started. For the boundary case of $|S|=n-1$, with just one task left being processed and no more tasks to schedule, we have that the worst-case remaining makespan is the maximum of all possible cases in which the newly scheduled task $k$ either completes before the currently processing task $i$, or after it: 
\begin{align}
 T'([S,k],F,D,i ,\bar{D}_{i}) = \max\left\{\max_{d_i:d\in U_{[S,k], F,D,i,\bar{D}_{i}}} d_i - \bar{D}_{i},\max_{d_k:d\in U_{[S,k], F,D,i,\bar{D}_{i}}} d_k\right\}. \label{eq:adversary_recursion_boundary}
\end{align}
An important feature of the adversary-perspective \gls{DP}, which is missing from the scheduler-perspective problem, is that its decision space is of finite dimension, and thus can be optimized over. From a mathematical optimization point of view, we can roughly say that the adversary \gls{DP} problem is a dual representation of the scheduler-perspective \gls{DP} problem, which, as it turns out, is easier to solve. In the following section, we present an \gls{MILO} reformulation of this \gls{DP}. %We placed the tedious extensions of the DP formulations into settings with more than two machines in Appendices \ref{app;DP_more_than_two} and \ref{app:DP_more_than_two_advers}.

\section{\revise{\gls{MILO} formulation and the \gls{2SSA} heuristic}} \label{sec:MIO}
\subsection{Mixed-integer problem formulation}
In this section, we formulate the adversarial perspective \gls{DP} of Section~\ref{sec:adversary.DP} as an \gls{MILO}. It is important to note that the adversary optimizes over an entire {\textit scenario tree} rather than on a {\it single} vector of task durations $d$ (\ie, a branch). This is because the adversary is also non-anticipative, in the sense that at each point in time, it makes a decision on the duration of the next task to complete while taking into account all possible future decisions of the scheduler.

To construct the scenario tree, we utilize two elements of the state description: $S$ and $F$ -- the ordered lists of tasks according to their starting and completion orders, respectively; we denote their combination by $\sigma=(S,F)$. We consider different states in which either (i) one of the tasks has just been completed and the scheduler needs to make a decision (in association with the inner minimum in the \gls{DP} recursion presented in  \eqref{eq:adversary_DP_m_machines}), or (ii) a new task has just been scheduled and the adversary decides which of the running tasks finishes first and what its remaining duration would be (associated with the outer maximization and inner maximizations in the \gls{DP} recursion presented in  \eqref{eq:adversary_DP_m_machines}). Before presenting the rigorous notation, let us discuss the main idea behind the scenario tree for a setting with four tasks.
\subsubsection{Illustrative case: Two machines and four tasks.} \label{sec:2machines4tasks}
Consider the following {\it states}:
\begin{itemize}
\item $\sigma = (S, F) = ([],[])$: Initial state, no tasks are scheduled yet -- both machines are idle; a task needs to be scheduled immediately.
\item $\sigma = ([1],[])$: Task 1 has been scheduled on one of the machines -- one machine is idle; a task needs to be scheduled.
\item $\sigma = ([1,2],[])$: Task 1 has been scheduled on one machine and task 2 on the other machine -- no machine is idle; the scheduler needs to wait until either task 1 or task 2 finish processing.
\item $\sigma = ([1,2],[2])$: Having initially scheduled tasks 1 and 2, task 2 finishes first -- a machine becomes idle; the scheduler has to schedule another task.
\item $\sigma = ([1,2,3],[2])$: Having initially scheduled tasks 1 and 2, task 2 finishes first and task 3 is scheduled immediately -- now the scheduler needs to wait until either task 1 or task 3 to finish processing.
\end{itemize}
The states evolve until both $S$ and $F$ contain all four tasks. To see how $\sigma=(S, F)$ defines the order of events, consider $\sigma=(S, F) =([1,2,3,4],[1,3,2,4])$, which describes an end state obtained after the following sequence of decisions/events:
\begin{itemize}
\item Tasks 1 and 2 are scheduled first,
\item task 1 finishes first,
\item task 3 is scheduled on the machine where task 1 was processing,
\item the next task to finish is task 3 (task 2 is still running),
\item Task 4 is scheduled on the newly available machine (where task 3 was processing), and
\item task 2 finishes followed by task 4, which finishes last.
\end{itemize}
Given the order of events as depicted by $\sigma$ (or more generally, a state within the scenario graph), we can formulate the inequalities that define the sequence of events by comparing the duration of the tasks scheduled on the different machines. Specifically, for the sequence presented above, the corresponding inequalities are:
\begin{align*}
d_1 & \leqslant d_2 \\
d_1 + d_3 & \leqslant d_2 \\
d_1 + d_3 + d_4 & \geqslant d_2,
\end{align*}
where the first inequality follows from the fact that tasks 1 and 2 are scheduled at the same moment but the duration of 1 is shorter, the second inequality follows from the fact that task 3 started immediately after task 1 but completed before task 2, and the third inequality follows from task 4 that started immediately after task 3 and completed after task 2. We can present the system of inequalities as:
$$
E_\sigma d_\sigma \leq 0, \quad E_\sigma = \left[ \begin{array}{rrrr}
1 & -1 & 0 & 0 \\
1 & -1 & 1 & 0 \\
-1 & 1 & -1 & -1 \\
\end{array}
\right], \quad d_\sigma= [d_1 \ d_2 \ d_3 \ d_4]^\top.
$$
Given $\sigma$, we can also express the total makespan in terms of the duration of different tasks. For example, denoting the makespan represented by $\sigma$ by $t_\sigma$, $t_\sigma = d_1 + d_3 + d_4$, or equivalently
$$
t_\sigma = e_\sigma^\top d_\sigma, \quad e_\sigma = \left[ \begin{array}{cccc} 1 & 0 & 1 & 1 \end{array} \right]^\top.
$$

\begin{figure}[H]
\centering
\begin{tikzpicture}[sibling distance=5em, level 1/.style ={level distance=1.3cm,sibling distance=4em}, level 2/.style ={level distance=1.3cm}, level 3/.style ={level distance=1.3cm}, level 4/.style ={level distance=2.5cm,sibling distance=7em},
  every node/.style = {shape=rectangle, draw, align=center, font=\scriptsize}]]
  \tikzset{greynode/.style={fill=black!30}}
  \node [greynode] {[],[]}
  	child { node {$[1,2],[]$}
      child { node [greynode] {$[1,2],[1]$} 
        child { node {$[1,2,3],[1]$} 
        	child {node {\scriptsize $[1,2,3,4], [1,3,2,4]$} edge from parent node[sloped, above, draw=none] {$d_3 + d_4$}}
        	child {node {$[1,2,3,4], [1,3,4,2]$} edge from parent node[sloped, above, draw=none] {$d_2 - d_1$}}
        	child {node {$[1,2,3,4], [1,2,3,4]$} edge from parent node[sloped, above, draw=none] {$d_2 + d_4 - d_1$}}
        	child {node {$[1,2,3,4], [1,2,4,3]$} edge from parent node[sloped, above, draw=none] {$d_3$}}} 
        child { node {[$1,2,4],[1]$} }  edge from parent node[sloped, above, draw=none,anchor=south,xshift=-3pt] {$d_1$} } 
      child { node [greynode] {$[1,2], [2]$} } }
    child {node {$[1,3],[]$}}
    child {node {$[1,4],[]$}}
    child {node {$[2,3],[]$}}
    child {node {$[2,4],[]$}}
    child {node {$[3,4],[]$}} ;
\end{tikzpicture}
\captionsetup{font=scriptsize}
\caption{A partial representation of a scenario tree for the case of $n = 4$ and $m = 2$. As the scheduling process progresses, one moves from the initial none $([],[])$ downwards. On gray nodes the scheduler decides to move to one of the direct children. The other nodes are adversary nodes, in which the adversary controls the branching direction. The labels on the edges starting from the adversary nodes denote the time it takes to move from one node to another.} \label{fig.states.tree}
\end{figure}
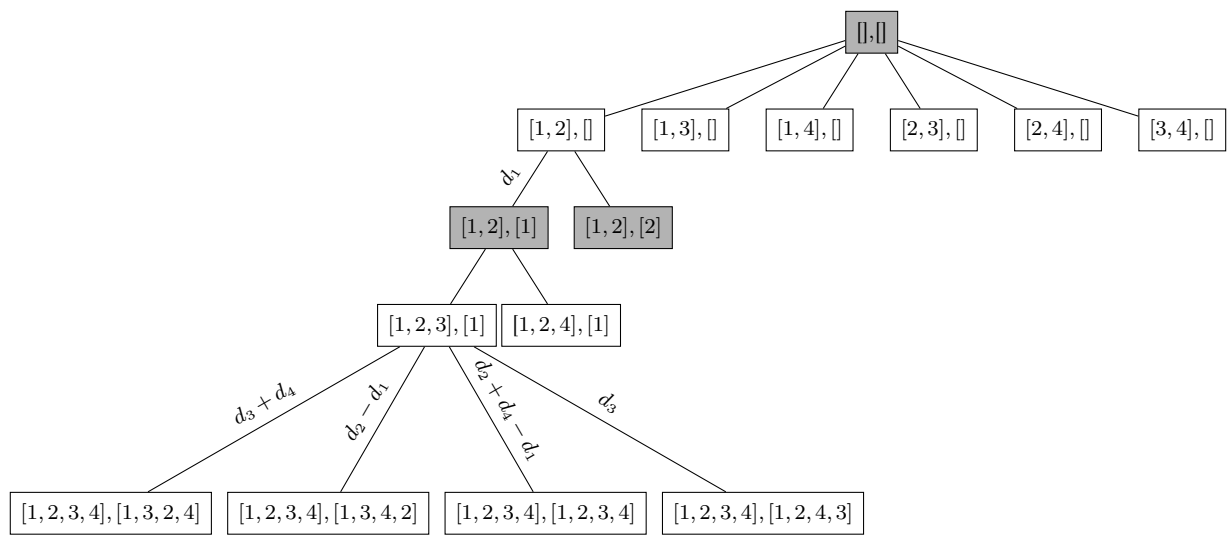

In a similar way, we can encode all the states within a scenario tree, where its root node represents the initial state $\sigma = ([],[])$ from which one proceeds to the next nodes by adding a task to either $S$ or $F$. Figure~\ref{fig.states.tree} illustrates a part of the tree. Note that we discarded states with $|S|=1$ when moving from the initial state since it is always optimal to schedule two tasks on the initially idle machines. In other words, state $([],[])$ (initial state) progresses directly to states in which $|S|=2$ such as $([1,2],[])$.
Similarly, we excluded states where there are no further decisions to be made by the scheduler or the adversary -- that is, states with $|S|=3$ and $|F|>1$ or $|S|=4$ and $|F|<4$. The last scheduler decision is the one after which a single task is left to be scheduled. After that, the adversary decides on the completion order of the tasks. In our setting of $n = 4$, the last state in which the scheduler makes a decision has two started tasks where one of them already completed. After this decision, three tasks have started and one completed and the adversary makes the final decisions about the completion times of the in-process tasks and the task not yet started (the last remaining task is automatically scheduled when a machine becomes available).

As shown in Figure~\ref{fig.states.tree}, task durations corresponding to certain root-to-leaf paths are `common'. For example, consider the two separate vectors of task duration corresponding to leaves $(S,F) = ([1,2,3,4],[1,3,2,4])$ and $(S,F) = ([1,2,3,4],[1,3,4,2])$. In both vectors, the  decision about the duration of task 1 was made before the last scheduling decision was made, and thus the duration of task 1 has to be the same for both paths.

Overall, in Figure~\ref{fig.states.tree}, the scheduler aims to progress from the root node to the lowest level in the tree such that
\begin{itemize}
\item at a gray node, the scheduler selects the direct children node to go to, and
\item at a white node, the adversary selects the direct children nodes to go to. 
\end{itemize}
The adversary's optimization problem is, therefore, to set the arc lengths in such a way that the shortest path for the scheduler is as long as possible. In the next section, we shall explicitly define this problem.

\subsubsection{Formulating the adversary's MILO problem.}
{We denote the set of all states $\sigma=(S,F)$} as $\mathcal{V}$. To decrease the state space of the scenario tree, certain states are discarded (see also the example described in Section~\ref{sec:2machines4tasks}). In particular,
\begin{itemize}
\item {We omit states where $|S| = 1$, since they correspond to an initialization of the system in which one machine is idle, thus an additional task must be scheduled immediately.} Hence, the states that follow the initial state are characterized by $|S|=2$ ({both} machines are busy).
\item Similarly, the last decision to be made is the one after which a single task remains to be scheduled. Such a decision is made at a state where the length of $S$ is $n - 2$ and the length of $F$ is $n - 3$. Afterwards, the adversary decides which tasks to complete and in which order. Hence, the states that follow states with  $|S|=n-1$ and $|F|=n - 3$ are the end states in which $|S|=|F|=n$.
\end{itemize}
Given this, we consider the state space $\mathcal{V}$ consisting of the sets of the scheduler states $\mathcal{D}$ and adversary states $\mathcal{N}$ (Table \ref{tab:scenario_states}).
\begin{center}
\begin{table}[ht]
\caption{The sets of states for the scheduler $\mathcal{D}$ (left) and the adversary $\mathcal{N}$ (right) when $m=2$.}\label{tab:scenario_states}
\vspace{10pt}
\begin{minipage}{.5\textwidth}
{\centering
{\def\arraystretch{1.2}\tabcolsep=3pt
{\scriptsize 
\begin{tabular}{p {0.45\textwidth} p {0.45\textwidth}} 
\hline
\textbf{$\boldsymbol{S}$-length} & \textbf{$\boldsymbol{F}$-length} \\ \hline
0 & 0 \\ \hline
$2$ & 1 \\\hline
$3$ & 2 \\\hline
\vdots & \vdots \\\hline
$n - 2$ & $n - 3$\\ \hline
\end{tabular}}}}\\
\end{minipage}
\begin{minipage}{.5\textwidth}
{\centering
{\def\arraystretch{1.2}\tabcolsep=3pt
{\scriptsize 
\begin{tabular}{p {0.45\textwidth} p {0.45\textwidth}} 
\hline
\textbf{$\boldsymbol{S}$-length} & \textbf{$\boldsymbol{F}$-length} \\ \hline
$2$ & 0 \\ \hline
$3$ & 1 \\\hline
\vdots & \vdots \\\hline
$n - 1$ & $n - 3$\\ \hline
$n$ & $n$  \\
\hline
\end{tabular}}}}\\
\end{minipage}
\end{table}
\end{center}
Let $t_\sigma$ denote the optimal worst-case makespan corresponding to being in state $\sigma \in \mathcal{D} \cup \mathcal{N}$. We denote the set of all final states of the tree, for which $|S|=|F|=n$, by $\mathcal{L}$. 
For each state $\sigma \in \mathcal{D}$, the scheduler selects the task that minimizes the worst-case remaining makespan (corresponding to the inner minimum of the \gls{DP} formulation in  \eqref{eq:adversary_DP_m_machines}). Thus, 
$$t_\sigma=\min\limits_{\delta \in \text{Children}(\sigma)} t_\delta,\; \sigma\in \mathcal{D}.$$
For each state $\sigma \in \mathcal{N}\setminus\mathcal{L}$, the adversary determines task completion times that maximize the worst-case remaining makespan (corresponding to the outer maximum within the \gls{DP} formulation in  \eqref{eq:adversary_DP_m_machines}). Thus, 
$$
t_\sigma=\max\limits_{\delta \in \text{Children}(\sigma)} t_\delta,\; \sigma\in \mathcal{N}\setminus\mathcal{L}.
$$
For each end state $\sigma \in \mathcal{L}$, we define a separate decision variable of the adversary being a vector $d_\sigma \in \mathbb{R}^n$. 

As illustrated in Section~\ref{sec:2machines4tasks}, for each final state $\sigma$, the set of inequalities that accommodates the task durations within the given scenario can be formulated in terms of the vector $d_\sigma$ as follows:
$$
E_\sigma d_\sigma \leq 0 , \ t_\sigma = e_\sigma^\top d_\sigma.
$$
To achieve non-anticipativity of the adversary decision with respect to future scheduler decisions, the adversary decision vectors $d_\sigma, d_\delta \in \mathbb{R}^n$ for $\sigma, \delta \in \mathcal{L}$ have to be the same as long as sequences $\sigma$ and $\delta$ are indistinguishable from each other in terms of scheduler decisions. Therefore, denoting $\mathcal{A}(\sigma,\delta)\in \mathcal{D}$ as the last common ancestor of both states $\sigma$ and $\delta$ before a different scheduling decision is made, we can define 
$$
d_\sigma(i)=d_\delta(i), \quad \forall i\in F, \quad (S,F)=\mathcal{A}(\sigma,\delta).
$$

The combination of task duration inequalities, $E_\delta d_\delta \leq 0$, and the non-anticipativity constraints can potentially lead to a situation where a realization chosen for some $\sigma\in \mathcal{L}$, denoted by $\bar{d}_\sigma$, makes the feasible set for another $\delta \in \mathcal{L}$ empty:
$$
\{ d \in U: \ E_\delta d_\delta \leq 0, \ d_\delta(i) = \bar{d}_\sigma(i), \quad \forall i\in F, \quad (S,F)=\mathcal{A}(\sigma,\delta) \} = \emptyset,
$$
\ie, a specific realization of the first (few) tasks makes ending up in a given $\delta \in \mathcal{L}$ impossible. We will use big-$M$ reformulations of $E_\sigma d_\sigma \leq 0$ to account for this possibility correctly.

We proceed to the problem formulation. If we denote by $t_0$ the optimal worst-case makespan corresponding to the root decision node $\sigma = ([],[])$, the adversary's problem \eqref{opt-P} can be stated as:
\begin{align} 
\max\limits_{d_\sigma, t_\sigma, z_\sigma} \ & t_0 \tag{P}\label{opt-P} \\
\text{s.t.} \ & t_\sigma = \min\limits_{\delta \in \text{Children}(\sigma)} t_\delta && \forall \sigma \in \mathcal{D} \label{eq:tree.dec.nodes} \\ 
& t_\sigma = \max\limits_{\delta \in \text{Children}(\sigma)} t_\delta && \forall \sigma \in \mathcal{N} \setminus \mathcal{L} \label{eq:tree.nat.nodes} \\
& d_\sigma(i)=d_\delta(i)&& \forall i\in F, (S,F)=\mathcal{A}(\sigma,\delta), \sigma,\delta \in \mathcal{L}  \label{eq:tree.nonanticipativity}\\%d_{\sigma, \sigma(i)}  = d_{\delta, \delta(i)} && \forall i \in I(\sigma, \delta), \forall \sigma,\delta \in \mathcal{L}  \label{eq:tree.nonanticipativity} \\
& t_\sigma = e_\sigma^\top d_\sigma - z_\sigma M && \forall \sigma \in \mathcal{L} \label{eq:tree.leaf.length} \\
& E_\sigma d_\sigma \leq z_\sigma M && \forall \sigma \in \mathcal{L} \label{eq:tree.leaf.set}\\
& z_\sigma\in\{0,1\}, \; d_\sigma \in U && \forall \sigma \in \mathcal{L},
\end{align}
where $M$ is a large positive number.
\begin{comment}
The constraints are:
\begin{itemize}
\item \eqref{eq:tree.dec.nodes}: The optimal makespan of the scheduler's node is the minimum of its direct children
\item \eqref{eq:tree.nat.nodes}: The optimal makespan of the adversary's node is the maximum of its direct children
\item \eqref{eq:tree.nonanticipativity}: The non-anticipativity constraints on the adversary.
\end{itemize}
\end{comment}
The big-$M$ constraints, \eqref{eq:tree.leaf.length} and  \eqref{eq:tree.leaf.set}, jointly with the binaries $z_\sigma$, enable one of two situations to occur for each of the corresponding paths $d_\sigma$ for $\sigma\in \mathcal{L}$: (i) $d_\sigma$ satisfies constraint \eqref{eq:tree.leaf.set} with $z_\sigma=0$ (which is possible depending on the uncertainty set structure and the non-anticipativity constraint \eqref{eq:tree.nonanticipativity}) and has a corresponding finite makespan, (ii) if constraint \eqref{eq:tree.leaf.set} cannot be satisfied, then $z_\sigma$ is set to $1$ and the makespan corresponding to this path will be $-\infty$, thus it will never qualify as the worst-case makespan. We also note that since the adversary solves a maximization problem, it will never allow all the children leaves of a certain scheduler decision branch to take a value of $-\infty$, since that would imply that this will be the value of the entire problem.

Because the problem is a maximization one, constraints \eqref{eq:tree.dec.nodes}, involving minimum terms, can be reformulated as
$$t_\sigma\leq  t_\delta,\quad \forall \delta \in \text{Children}(\sigma),\quad \sigma\in\mathcal{D},$$
while constraints \eqref{eq:tree.nat.nodes}, involving maximum terms, need to be implemented using auxiliary binary variables and a big-$M$ constraint as follows:
\begin{align}
t_\sigma&\leq t_\delta+M\cdot w_{\sigma,\delta},&& \forall \delta \in \text{Children}(\sigma),\sigma\in\mathcal{N}\setminus\mathcal{L}\\
\sum_{\delta\in \text{Children}(\sigma)} w_{\sigma,\delta}& \leq |\text{Children}(\sigma)| - 1 && \forall \sigma\in\mathcal{N}\setminus\mathcal{L}\\
w_{\sigma,\delta}&\in\{0,1\},&& \forall \delta \in \text{Children}(\sigma),\sigma\in\mathcal{N}\setminus\mathcal{L}.
\end{align}
%\subsubsection{Complexity of the \texorpdfstring{\gls{MILO}} \ \ formulation.}
In Appendix \ref{More_than_two_MIO}, we present the \gls{MILO} formulation for settings with more than two machines. Although exact, this formulation's complexity is exponential, as shown in Appendix~\ref{appendix:milo.complexity}. For that reason, in the next section we develop a scalable adaptive heuristic. I
\subsection{\revise{Two-stage \gls{SA} heuristic (\gls{2SSA})}}
\revise{As the exact \gls{MILO} formulation of \gls{AR} does not scale well with the problem size, we propose a scalable heuristic for the two machine setting ($m=2$) that takes into account only one adaption of the task allocations, as opposed to $n - 2$ such adaptations in the exact reformulation.}

\revise{The idea of the heuristic is as follows: we schedule two tasks in parallel initially and, based on which task finishes first and its duration, apply the \gls{SA} to the remaining tasks. The schedule of the remaining tasks adapts itself thus only to the durations of the two initial tasks, remaining fixed afterwards.} 

\revise{The heuristic divides the time into two parts: before and after one of the two initial tasks completes; intuitively, the heuristic strives to place tasks whose durations are most informative about the remaining task durations as the initial tasks.}

\revise{For the heuristic we need to solve a two-stage robust optimization problem. To formulate it, we define $x \in \{0, 1 \}^n$ to be a vector of task allocations as in Section~\ref{sec:bound.promised.durations}. The two-stage problem is given by:
\begin{align} \label{eq:two.stage.SA}
\min\limits_{i,j \in[n]: i < j} \max \left\{ \begin{array}{l} \sup\limits_{\tilde{d}_i: \ \tilde{d} \in U,\  \tilde{d}_i < \tilde{d}_j} \pinf\limits_{\substack{x \in \{0, 1 \}^n: \\ x_i = 1, \ x_j = 0}} \sup\limits_{\substack{d \in U: \ d_i=\tilde{d}_i,\\d_i < d_j}} x^\top d, \\ \sup\limits_{\tilde{d}_j: \ \tilde{d} \in U,\  \tilde{d}_i \geq \tilde{d}_j} \pinf\limits_{\substack{x \in \{0, 1 \}^n: \\x_i = 1, \ x_j = 0}} \sup\limits_{\substack{d \in U:\  d_j=\tilde{d}_j\\ d_i \geq d_j}} (e - x)^\top d \end{array}  \right\}.
\end{align}
The first minimum determines which two tasks $i$ and $j$ are allocated first by the scheduler. Without loss of generality, we assume task $i$ is scheduled on the first machine and task $j$ on the second machine. Next, via the maximum and the first supremum, the adversary decides which of these two tasks finishes first and with what duration. Based on this information, in the next minimum, the scheduler decides on the \gls{SA} schedule for the remaining tasks and the adversary decides (the last supremum) on the exact durations of the remaining tasks.}

\revise{For fixed $i$, $j$, this problem can be solved in a tractable way using the Column-and-Constraint generation algorithm of \cite{zeng2013solving} for two-stage robust optimization problems, for $n\leq 20$. In Appendix~\ref{app:ccg}, we outline the actual \gls{2SSA} implementation; to apply the heuristic one needs to solve $n(n-1)/2$ smaller problems (\ie, one for each $i,j\in[n]$ such that $i< j$).}

\section{Numerical study} \label{sec:numerical.study}
Our analytical work is accompanied by a numerical study that focuses on problem instances with two machines. We start with five-task instances, which are large enough to differentiate between alternative policies and elicit some insights, yet small enough for to allow us to run the \gls{MILO} formulation presented in the last section. We use these instances to demonstrate the value of adaptivity, and to compare the \gls{AR} and \gls{2SSA} performance. We then proceed to investigate larger instances with 10-20 tasks, comparing the performance of \gls{SA} to that of \gls{2SSA} in favor of validating our theoretical bounds in terms of identifying the settings in which accounting for adaptivity in the planning stage leads to a significant benefit.
 
\subsection{The setup} 
We investigate the  \gls{SA}, \gls{SL}, \gls{AR} and \gls{2SSA} policies. For a realistic comparison, we implement the policies for each scenario in a rolling-horizon fashion similarly to the way they are expected to be applied in reality. \revise{We note that while we implemented the policies without the re-optimizing via rolling-horizon we do not include the related results, since they were consistently inferior compared to the rolling-horizon implementations.} A scenario represents a particular realization of task durations for the tasks. We conduct each experiment as follows: At the beginning of the planning horizon, when task durations are only known to belong to $U$, we select the optimal two first tasks to be scheduled, according to the considered policy. If multiple initial scheduling decisions give the same worst-case makespan, we choose between them according to lexicographic ordering (w.r.t. task indices). Then, we determine which of the two running tasks finishes first (based on the known duration for the considered scenario). Next, we update the uncertainty set to include the information about the duration of the completed task and for how long the other task has been processing so far, and calculate the optimal next scheduling decisions by the considered policy. For \gls{AR}, for example, this means solving the \gls{MILO} using the current information about task duration and completion times in order to select the task to be scheduled for the idle machine. Once both machines are busy, we again determine which of the running tasks finishes first, and optimally select the next task to be scheduled for the first-to-be-idle machine. This sequence of events continues until all tasks have been allocated and finished.

We find a lower bound on the obtainable makespan, by applying a \gls{PH} policy for each scenario. Namely, we assume that the specific scenario is known in advance and determine the optimal task-to-machine allocations that minimize the makespan. The makespan obtained by the \gls{PH} policy is a lower bound on the possible makespan achieved by any policy. Therefore, we use the \gls{PH} makespan as a benchmark for the investigated policies.

To define the performance measures, we denote $p\in\{\gls{AR},\gls{SL},\gls{SA},\gls{2SSA},\gls{PH}\}$ as the selected policy, $k\in\{1,\dots,N\}$ identifies a problem instance and $s\in\{1,\dots,R\}$ denotes the scenario number. Accordingly, we denote $\tilde{T}_{k,p,\cdot}$ and  $T_{k,p,s}$ for problem instance $k$, policy $p$, and scenario $s$, as the promised worst-case makespan at the beginning of the planning horizon and the realized makespan when the policy is applied via a rolling-horizon approach, respectively. Notice that the promised worst-case makespan  $\tilde{T}_{k,p,\cdot}$ is independent of the scenario's realization and is only affected by the specification of the uncertainty set, since it is the makespan for the case in which the worst-case scenario is realized.  

The mathematical formulas and names of the measures that we report are described in Table~\ref{tab:results.measures}. % and the results are given in Tables~\ref{tab:results.box} and \ref{tab:results.ellipsoid}. 
 \revise{The first measure is the percentage of instances where the policy's initial decision was not optimal for the \gls{AR} policy.}
The next three measures specify the actual values (max and average) of the makespan, for comparison purposes. The final four measures give us an indication of how well the policies perform relative to the \gls{AR} and \gls{PH} solutions.
\begin{itemize}
\item {\bf Promised to max \gls{PH}.} The upper bound on the ``price'' that a risk-averse scheduler who commits (\eg,  in a contract) to the promised makespan is expected to pay with respect to a \gls{PH} policy. The reasoning behind this measure is that before any of the tasks starts processing, the scheduler does not know which scenario will be realized, thus she compares the promised makespan to the maximal perfect hindsight makespan across scenarios. 
\item {\bf Max makespan to max \gls{PH}.} Following the previous measure, a ratio between the maximal realized makespan to the maximal perfect hindsight makespan gives an indication into the expected price that a risk averse scheduler who applies the rolling-horizon policies would contractually pay with respect to the \gls{PH} policy.
\item {\bf Makespan to \gls{PH}.} While \gls{RO} does not optimize for non-worst-case scenarios, we approximate the average ratio for each scenario between the realized makespan by the applied policy and its perfect hindsight counterpart. This measure gives an indication into the expected price of the risk-averse scheduler for an ``average'' scenario when compared with the corresponding perfect hindsight makespan.
\item {\bf Max makespan to max \gls{AR}.} To benchmark the best robust policy, which is \gls{AR} with its alternatives -- \gls{SA}, \gls{SL} and \gls{2SSA} -- we replicate the ``Max makespan to max \gls{PH}'' measure with \gls{AR} replacing \gls{PH} as the reference point.
%\item {\bf Max of max makespan to max \gls{AR}.} In order to observe the maximal possible difference between the policies and the \gls{AR} policy, we also present the maximal ratio over all random instances of the worst-case realized makespans.
\end{itemize}
We note that these experiments evaluate the performance of the policies under conditions of worst-case (denoted as the promised makespan) and non-worst-case scenarios in which we uniformly sample task durations from the respective scenario uncertainty sets. It is important to experiment with both conditions since the \gls{RO} methodology is indifferent with respect to non-worst-case scenarios, thus an optimal robust policy may perform poorly when applied on such scenarios \citep[see,][]{iancu2014pareto}.
\begin{table}[ht!]
\centering
\caption{Description and formulas of the performance measures} \label{tab:results.measures}
{\def\arraystretch{1.2}\tabcolsep=3pt
\scriptsize {\begin{tabular}{lll}\hline
\textbf{Short name} & \textbf{Description} & \textbf{Calculation}\\ \hline
\revise{Suboptimal initial decision} & \revise{Percentage of instances with initial decisions different from \gls{AR}} & - \\ \hline
Promised makespan & Average promised makespan of policy $p$& $\frac{1}{N}\sum_{k=1}^{N} \tilde{T}_{k,p,\cdot}$ \\  \hline
Max makespan&Average maximal makespan of policy $p$&  $\frac{1}{N}\sum_{k=1}^{N}\max_{s\in [R]} T_{k,p,s}$ \\  \hline
%Max \gls{PH} makespan&Average maximal makespan under perfect hindsight & $\frac{1}{N}\sum_{k=1}^{N}\max_{s\in [R]} T_{k,\gls{PH},s}$ \\  \hline
Makespan&Average makespan of policy $p$& $\frac{1}{N}\sum_{k=1}^{N} \frac{1}{R}\sum_{s=1}^R T_{k,p,s}$ \\  \hline
%\gls{PH} makespan&Average perfect hindsight makespan &  $\frac{1}{N}\sum_{k=1}^{N} \frac{1}{R}\sum_{s=1}^R T_{k,\gls{PH},s}$ \\  \hline
Promised to max \gls{PH}&\makecell[l]{Ratio between the promised of policy $p$ \\ to the maximal makespan under perfect hindsight}& %$\frac{1}{N}\sum_{k=1}^{N}
$\frac{\tilde{T}_{k,p,\cdot}}{\max_{s\in[R]} T_{k,\gls{PH},s}}-1$\\  \hline
Max makespan to max \gls{PH}&\makecell[l]{Ratio between the maximal of policy $p$ \\ to the maximal makespan under perfect hindsight}& %$\frac{1}{N}\sum_{k=1}^{N}
$\frac{\max_{s\in[R]} T_{k,p,s}}{\max_{s\in[R]} T_{k,\gls{PH},s}}-1$\\  \hline
Makespan to \gls{PH}&\makecell[l]{Ratio between the makespan of policy $p$\\ to the makespan under perfect hindsight} &  %$\frac{1}{N}\sum_{k=1}^{N}
$\frac{1}{R}\sum_{s=1}^{R}\frac{T_{k,p,s}}{T_{k,\gls{PH},s}}-1$\\  \hline
Max makespan to max \gls{AR}&\makecell[l]{Ratio between the maximal makespan of policy $p$\\ to the maximal makespan the \gls{AR} policy}& %$\frac{1}{N}\sum_{k=1}^{N}
$\frac{\max_{s\in [R]} T_{k,p,s}}{\max_{s\in[R]} T_{k,\gls{AR},s}}-1$\\  \hline
%Max of max makespan to max \gls{AR}&\makecell[l]{Maximal ratio between the maximal makespan of policy $p$\\ to the maximal makespan of the \gls{AR} policy} & $\max_{k\in[N]}\frac{\max_{s\in[R]} T_{i,p,s}}{\max_{s\in[R]} T_{k,\gls{AR},s}}$\\ \hline
\end{tabular}}}
\end{table}

\revise{The code was run on a PowerEdge R740xd server with two Intel Xeon Cold 6254 3.1GHz processors, each with 18 cores, and a total RAM of 384GB.}

\subsection{Small problem instances}
First, we test the policies over $N=500$ problem instances, where each instance $k$ has a discrete uncertainty set $U^k \subseteq \mathbb{R}^n$ ($n=5$) consisting of $|U^k|=R = 15$ scenarios. Each problem instance $k\in\{1,\dots,N\}$ is generated by drawing a vector of the nominal durations $d^{0,k}$ and their perturbations $\bar{d}^k$ out of a uniform distribution according to  $d^{0,k}=(d^{0,k}_1,d^{0,k}_2,\dots,d^{0,k}_n)\sim \text{Unif}([0.1, 2.0]^n)$ and  $\bar{d}^k=(\bar{d}^k_1,\bar{d}^k_2,\dots,\bar{d}^k_n)\sim \text{Unif}([0.1, 5.0]^n)$, respectively. Each of the $s\in\{1,\dots,R\}$ scenarios of instance $k$ is sampled as
$$
\tilde{d}_j^{k,s} = d^{0,k}_j + u_j^{k,s} \bar{d}^k_j, \ \ j = 1,\ldots, n,
$$
with $u^{k,s} \in \mathbb{R}^n$ sampled uniformly from one of the following sets:
\begin{enumerate}[(\Roman*)]
\item An $n$-dimensional ball $\{ u \in \mathbb{R}^n_{+}: \ \| u \|_2 \leq 1 \}$ (type I).
    \item An $n$-dimensional box $\{ u: \ \| u \|_\infty \leq 1 \}$ (type II)
\end{enumerate} 
We rounded $\tilde{d}^{k,s}$ to multiples of 0.1 to allow for the situation that two tasks finish simultaneously. We refer to the resulting discrete uncertainty sets as uncertainty sets of type I and II based on the method in which we generate $u^{k,s}$. Because the results for both set types are very similar, we discuss here only the type I sets, and relegating the type II results in Appendix~\ref{app:box_uncertainty_sampled}.

In Table~\ref{tab:results.ellipsoid} we summarise the first 4 measures comparing the different methods. The significance of \gls{AR} is underlined by the fact that its first-stage decisions were different from the other policies in a substantial portion of the problems. %For \gls{SA}, \gls{SL}, and \gls{2SSA} the percentages of different first-stage decisions, compared to \gls{AR}, are $50\%$, $11\%$ and $7\%$, respectively. 
Indeed, we found that there is a high correlation between percent of different initial decisions compared to \gls{AR} and the policies performances compared to \gls{AR}. That is, \gls{2SSA} is expected to be the best performing (and scalable) policy,  followed closely by \gls{SA} and \gls{SA} is last.

% \begin{table}[ht!]
% \TABLE
% {Performance measure results for type II uncertainty sets. \label{tab:results.ellipsoid}}
% {\centering
% {\def\arraystretch{1.2}\tabcolsep=3pt
% {\scriptsize \begin{tabular}{lccc} \hline
% \textbf{Short name} & \gls{AR} & \gls{SL} & \gls{SA}  \\  \hline
% Promised makespan & $ \mathbf{5.661} (0.927) $ & $ 5.700 (0.941) $ & $ 6.210 (1.051) $  \\  \hline
% Max makespan &  $ \mathbf{5.661} (0.927) $ & $ 5.681 (0.930) $ & $ 5.855 (0.976) $ \\  \hline
% Max \gls{PH} makespan & $ 5.656 (0.925) $ & $ 5.656 (0.925) $ & $ 5.656 (0.925) $ \\ \hline
% Makespan & $ 4.872 (1.010) $ & $ \mathbf{4.870} (1.011) $ & $ 4.913 (1.036) $  \\ \hline
% \gls{PH} makespan& $ 4.783 (1.008) $ & $ 4.783 (1.008) $ & $ 4.783 (1.008) $ \\ \hline
% \makecell[l]{Promised to max \gls{PH} (\%)} & $\mathbf{0.1} (  0.5) $ & $0.8 (1.5) $ & $9.8 (  5.8) $  \\  \hline
% \makecell[l]{Max makespan to max \gls{PH} (\%)} & $\mathbf{0.1} (0.5) $ & $0.5 (1.2) $ & $3.6 (4.5) $   \\  \hline
% \makecell[l]{Makespan to \gls{PH} (\%)} & $\mathbf{2.0} (3.9) $ & $\mathbf{2.0} (3.8) $ & $ 2.8 (4.9) $ \\  \hline
% \makecell[l]{Max makespan to max \gls{AR} (\%)} & $- $ & $0.4 (  1.1) $ & $3.5 (  4.4) $   \\  \hline
% \makecell[l]{Max of max makespan to max \gls{AR} (\%)}  & $- $ & $7.5 (-) $ & $31.4 (-) $\\ \hline
% \end{tabular}}}}{Makespan values are presented with their standard deviations in parentheses. Ratio values represent the difference, in percentage, with respect to the value of the policy in the denominator of the ratio and their standard deviations in parentheses.}
% \end{table}

\begin{table}[ht!]
\TABLE
{Performance measure results for type I uncertainty sets. \label{tab:results.ellipsoid}}
{\centering
{\def\arraystretch{1.2}\tabcolsep=3pt
{\scriptsize \begin{tabular}{lccccc} \hline
\textbf{Short name} & \gls{AR} & \gls{2SSA} & \gls{SL} & \gls{SA}  &\gls{PH}\\  \hline
\revise{Suboptimal initial} & - & \revise{$7\%$} & \revise{$11\%$} & \revise{$50\%$} & - \\ \hline
Promised makespan & $ 5.628 (0.914) $ & $ 5.653 (0.913) $ &  $ 5.668 (0.928) $ & $ 6.174 (1.034) $ &  - \\ \hline 
Max makespan & $ 5.628 (0.914) $ & $ 5.641 (0.913) $ & $ 5.648 (0.918) $ & $ 5.824 (0.962) $ & $5.624 (0.913)$ \\ \hline 
Makespan & $ 4.832 (1.013) $ & $ 4.821 (1.015) $ & $ 4.830 (1.017) $ & $ 4.874 (1.042) $  &  $4.743 (1.015)$\\ \hline 
\end{tabular}}}}{Makespan values are presented with their standard deviations in parentheses. }%Ratio values represent the difference, in percentage, with respect to the value of the policy in the denominator of the ratio and their standard deviations in parentheses.}
\end{table}

Our results indicate that the \gls{AR} promised makespan was the shortest, with \gls{2SSA} nearly the same (longer by $0.4\%$), \gls{SL} very close behind (longer by $0.7\%$) and \gls{SA} trailing behind with a longer makespan by $9.7\%$. The upper bound on the ``price'' that a risk-averse scheduler, who commits to the promised makespan, is expected to pay upfront with respect to a \gls{PH} policy (as indicated by the Promised to max \gls{PH} measure) was very small for \gls{AR} ($0.1\%$), \gls{2SSA} ($0.5\%$), \gls{SL} ($0.8\%$), and much higher for \gls{SA} ($3.6\%$).

%The makespan to \gls{PH} measure reveals that the difference between the robust makespans and their perfect hindsight counterparts is between $2-3.2\%$ for both types of uncertainty sets. This small difference with respect to a perfect hindsight policy indicates that robust policies for our \gls{PMS} instances are surprisingly good on average. 

%The max of max makespan to max \gls{AR} measure provides an indication of the highest difference that we can expect (for an unknown instance) between the applied policy and \gls{AR}. The makespans by \gls{SL} were higher by $4.9\%$ and $7.5\%$ for uncertainty type I and II sets, respectively. The makespans of \gls{SA} with respect to \gls{AR} were the highest---$17.4\%$ and $31.4\%$ for uncertainty type I and II sets, respectively. Thus, showing a possible significant advantage of using \gls{AR} over the other methods, even when applying them in a rolling-horizon fashion. 

\begin{figure}
\centering
\includegraphics[width=0.9\textwidth]{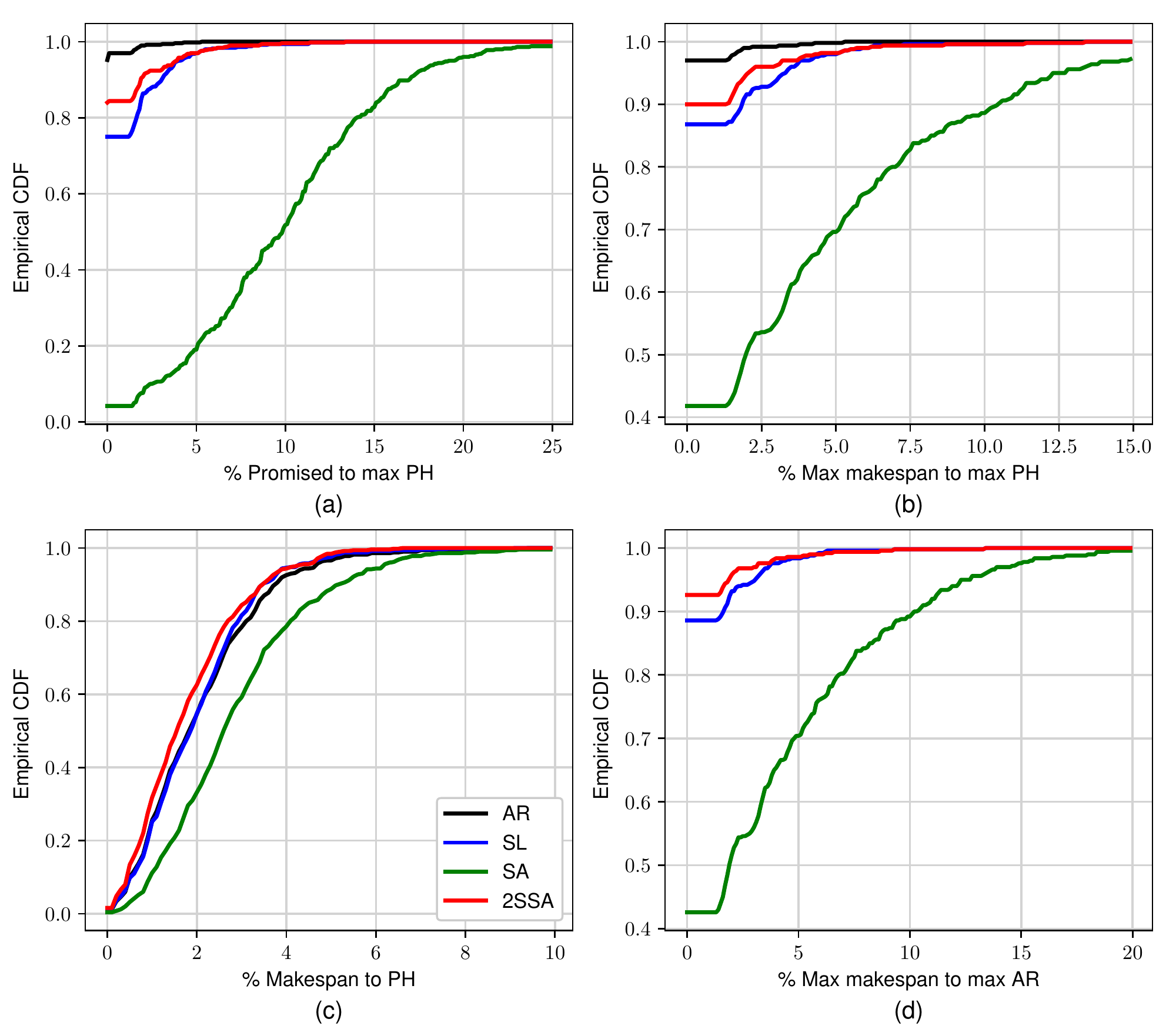}
\caption{\revise{Empirical \glspl{CDF} (over instances) of the last four performance measures of Table~\ref{tab:results.measures} for Type I uncertainty set.}}
\label{fig:ecdfs_type_1}
\end{figure}

Figure~\ref{fig:ecdfs_type_1} presents the \gls{CDF} of the last four performance measures in Table~\ref{tab:results.measures}, by means of empirical \glspl{CDF} across instances. Figure~\ref{fig:ecdfs_type_1}(a) shows the inferior performance of \gls{SA} compared to the other policies in providing good promised makespans. Among the remaining policies, \gls{AR} is the best, followed by \gls{2SSA} and \gls{SL}.
Figure~\ref{fig:ecdfs_type_1}(b) shows similar relationships between the different methods with respect to the ratio of the worst-realization of a given policy compared to the worst \gls{PH} duration, with the \gls{SA} being clearly the worst.
When taking the less conservative measure of taking the average ratio of scenarios (Figure~\ref{fig:ecdfs_type_1}(c)), the four policies are closer to each other, with \gls{2SSA} becoming surprisingly the best.
In the end, in Figure~\ref{fig:ecdfs_type_1}(d) we see the worst-case performance of the policies \gls{SL}, \gls{SA} and \gls{2SSA} compared to \gls{AR}. While \gls{2SSA} and \gls{SL} with comparable performance to that of the \gls{AR}, the \gls{SA} can go as much as 20\% worse. %We expect that the relatively good performance of \gls{SL} and \gls{2SSA} is due to the small problem size.

% \begin{figure}
% \centering
% \includegraphics[width=0.7\textwidth]{plot2.pdf}
% \caption{Histograms (over instances) of the average percent the  makespan of \gls{AR}, \gls{SL}, \gls{SA} exceeds that of \gls{PH} (associated with makespan to \gls{PH}). }
% \label{fig:histograms_to_PH}
% \end{figure}

%Figure~\ref{fig:histograms_to_AR} illustrates the optimality loss of \gls{SL} and \gls{SA}, compared to \gls{AR} due to lack of full adjustability in the scheduling algorithm. From here, it is clear that \gls{SA} is substantially inferior to the \gls{AR} and \gls{SL}, as seen from the long tail of the distribution of the makespan differences. As mentioned earlier, the performance of \gls{AR} and \gls{SL} are comparable for both uncertainty sets; yet, by observing the \gls{SL} tail, we can deduce that  when ``bad'' scenarios are realized, \gls{AR} may have a significant advantage over its alternatives.

% \begin{figure}
% \centering
% \includegraphics[width=0.7\textwidth]{plot1.pdf}
% \caption{Histograms (over instances) of the percent the worst-case makespan of \gls{SL} and \gls{SA} exceeds that of \gls{AR} (associated with max makespan to max \gls{AR}).}
% \label{fig:histograms_to_AR}
% \end{figure}

To recapitulate, \gls{AR} is the best policy for the risk-averse scheduler but it is not scalable. Our scalable heuristic, the \gls{2SSA}, and the non-scalable \gls{SL} achieve comparable performance with respect to  \gls{AR}, while \gls{SA} demonstrates the worse performance. as expected, the gaps between all policies are smaller when applying the policies on average (non-worst-case) scenarios, yet \gls{SA} is still inferior with respect to the other policies. 

\subsection{\revise{Large problem instances}}
\revise{In this set of experiments, we test the policies for $n=10,15,20$ tasks. Because of the problem sizes we could only test the \gls{SA} and \gls{2SSA} policies. For each value of $n$, and budget $\Gamma=0.1n,0.2n,0.3n,0.4n,0.5n,0.6n$, we generated $N=50$ problem instances. For each instance, we used a budgeted uncertainty set $U^k \subseteq \mathbb{R}^n$ of the form
$$U^k=\left\{d\in\mathbb{R}^n: d_i=d_i^{0,k}+u_i\bar{d}_i^k,\ i\in[n],\ u\in[0,1]^n, \sum_{i=1}^n u_i\leq\Gamma\right\},$$
where the nominal durations $d^{0,k}$ and their perturbations $\bar{d}^k$ are randomly generated such that $d_i^{0,k}\sim\text{Unif}(0.5, 5.0)$ and  $\bar{d}^k_i\sim\text{Unif}(0.5, 1.0) \cdot d_i^0$, respectively for each $i\in[n]$. For each realized scenario, we found the \gls{PH} solution by assuming that the realized durations are known in advance.}

\begin{figure}[t]
    \centering
    \includegraphics[scale=0.7]{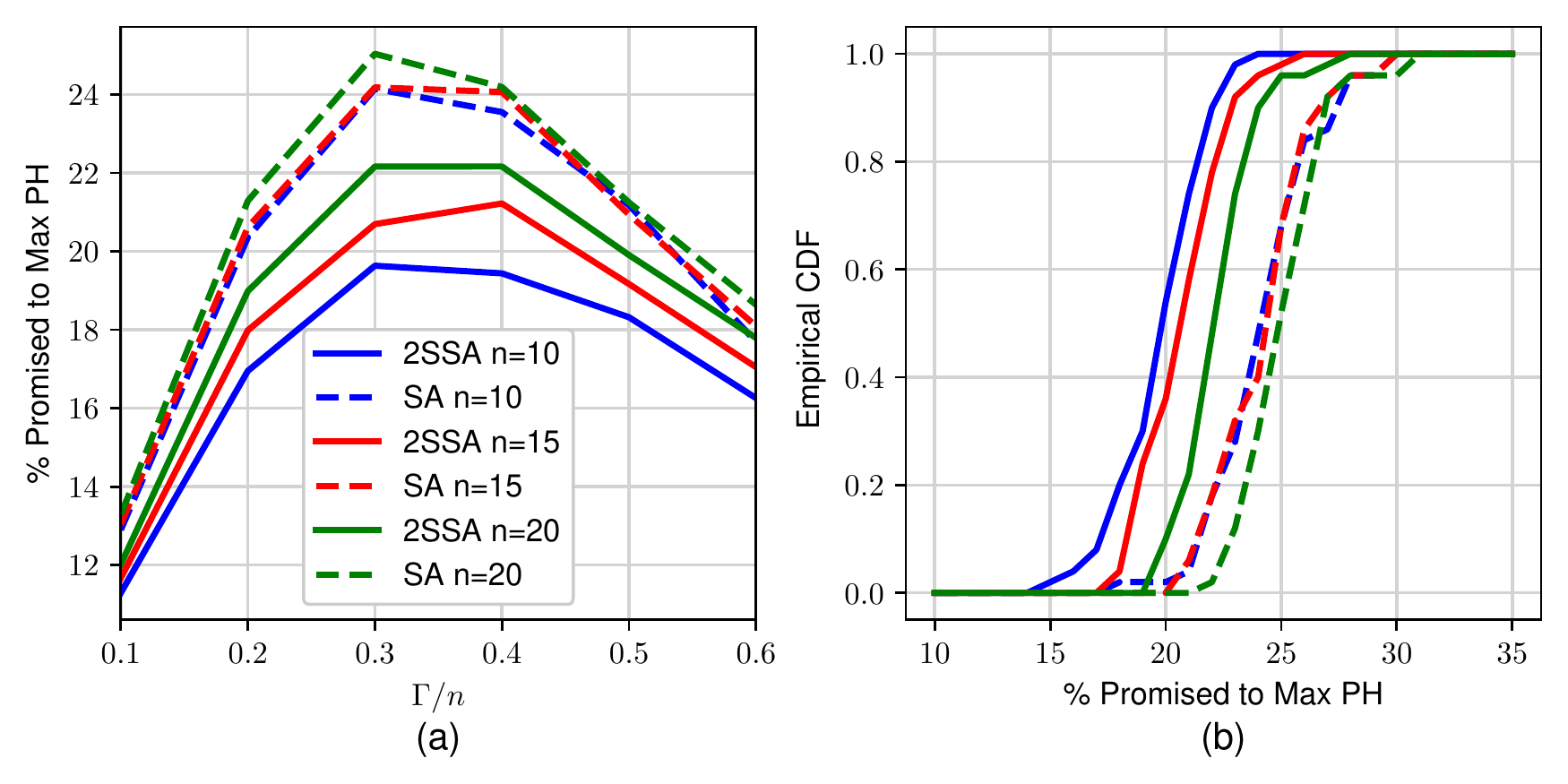}
    \caption{\revise{Promised to Max PH ratio for $n=10,15,20$ over $N=50$ instances: (a) average for different values of $\Gamma$ (b) empirical \gls{CDF} for $\Gamma=0.3n$.}}
    \label{fig:large_promised}
\end{figure}

\revise{We start by comparing the promised duration of \gls{SA} and \gls{2SSA} to the worst case duration of \gls{PH} on a sample of $K=50$ scenarios, which we present in Figure~\ref{fig:large_promised}. We observe the same phenomenon predicted by our bound in Section~\ref{sec:bound.promised.durations}, \ie, \gls{SA} promised duration gets close to those of the \gls{PH} when $\Gamma/n$ gets close to $0$ or $1$, and is further away from \gls{PH} when $\Gamma/n\in[0.3,0.4]$ in Figure~\ref{fig:large_promised}(c). We also note that this maximal difference is around $22\%$ implying a potential for improvement by methods which take into account adaptivity. Indeed, we see that the promised makespan of the \gls{2SSA} heuristic improves upon \gls{SA} in this region, by $2\%-4\%$. Moreover, as $n$ increases both methods' have worse promised duration with respect to the \gls{PH}, and the gap between them decreases. We explain the latter by the limited way in which \gls{2SSA} accounts for adaptivity in the planning stage. Since only one stage of adaptivity is taken into account, as $n$ increases there are more stages for which we do not account for adaptivity and thus, the relative benefit decreases. We note that the increase in the promised to \gls{PH} ratio for both methods as $n$ grows larger, is in contrast to their actual performance when implemented in a \gls{RH} fashion, as we discuss next. This gap between the promised and actual performance due to not taking into account enough of the future adaptivity may put a decision maker at a disadvantage when contractually committing to a specific makespan.}

\begin{figure}
    \centering
    \includegraphics[width=\textwidth]{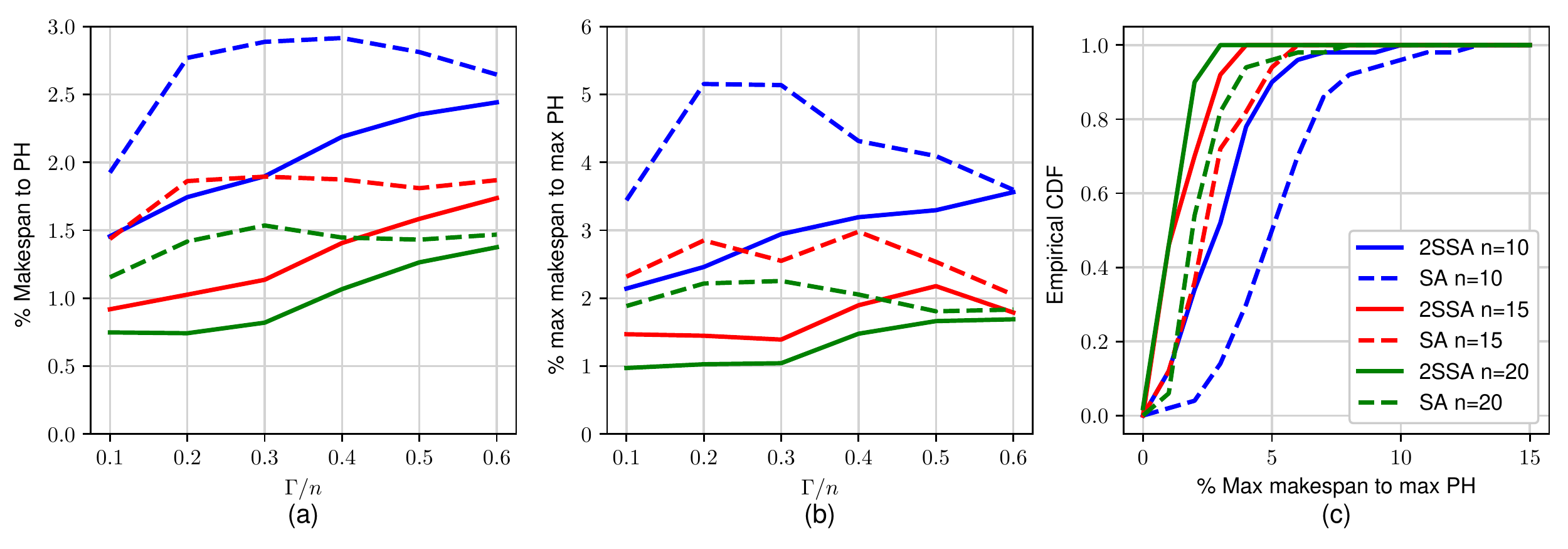}
    \caption{\revise{Makespan to PH ratios for $n=10,15,20$ over $N=50$ instances: (a) average `Makespan to PH' (b) average `Max makespan to max PH' (c) \gls{CDF} of `Max makespan to max PH' for $\Gamma=0.3n$.}}
    \label{fig:large_actual}
\end{figure}

\revise{To observe the actual performance of the methods, as they are applied in a \gls{RH} fashion, we generated $K=50$ random scenarios and computed, for each, the makespan to \gls{PH} ratio. Figure~\ref{fig:large_actual} depicts the performance measures of \gls{RH}-\gls{SA} and \gls{RH}-\gls{2SSA} as a function of both $n$ and $\Gamma$. Additionally, we present the empirical \gls{CDF} of `Max makespan to max PH' for $\Gamma=0.3n$. We see that as $n$ increases both the average and maximum gap between \gls{RH}-\gls{SA} and \gls{PH} decreases, as predicted by the bound in Section~\ref{sec:bound.rh}. In Figure~\ref{fig:large_actual}(b) we observe that for the lower values of $\Gamma/n$, \gls{2SSA} reduces the suboptimality of \gls{SA} w.r.t. the \gls{PH} by roughly 50\%.}

\revise{Importantly, Figure~\ref{fig:large_actual} demonstrates that, even in the rolling horizon implementation, \gls{2SSA} has favorable performance compared to the \gls{SA} for all tested settings.}

\revise{\subsection{Computational times}\label{sec:computational_times}}

\revise{We compared the computational times of our two scalable policies, the SA and 2SSA over instances with $n=5, 10, 15, 20$ tasks and under budgeted uncertainty. As the 2SSA can be easily parallelized over the $n(n-1)/2$ pairs of tasks to be scheduled first, we tested the running times with and without 2SSA parallelization. Figure~\ref{fig:times} summarizes the average run times of the entire rolling-horizon simulations for different $\Gamma / n$ ratios. It indicates that, given parallelization, the performances of SA and 2SSA are comparable for the tested problem sizes. We note that the presented running times for both \gls{SA} and \gls{2SSA} include all the re-optimization rounds over the entire rolling horizon. Thus, for the case of $n$ tasks, the time to adapt the scheduling decision each time a task completes is, roughly, the presented time divided by the $n-2$ times that each policy is re-optimized. These times are expected to be reasonable for most real-world settings, in the sense that machine will be idle only a short time while the problem is resolved. Practically, even these idle times can be decreased in classical manufacturing processes; this is due to the fact that the machines can accurately predict residual processing times when a part is close to completing its processing, allowing to solve the next re-optimization ahead of time.}   

\revise{\begin{figure}[ht]
    \centering
     \includegraphics[width=\textwidth]{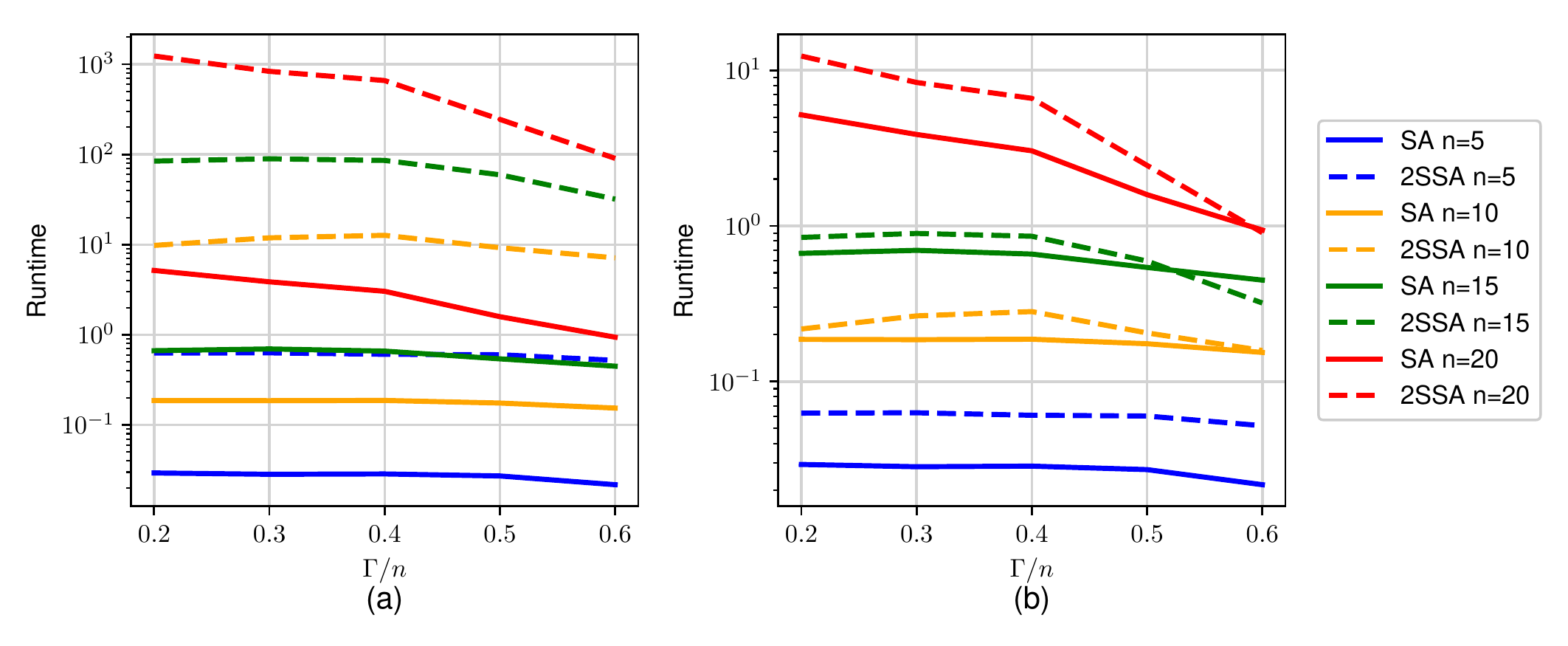}
    \caption{\revise{Computational times for  $n=5,10,15,20$ over the entire rolling horizon: (a) average `raw' computational times (b) average computational times with parallelization on up to 100 cores.}}
    \label{fig:times}
\end{figure}}
\section{Managerial insights}\label{sec:insights}
We outline our main managerial insights for the studied setting. The insights are relevant to schedulers within multiple domains that can be modeled via \gls{PMS} such as production lines in which machines process a set of tasks, computer multiprocessors (``cloud processing'') for processing jobs, shipyards and ports in which ships are loaded and unloaded, doctors who treat patients in a walk-in clinic or triage setting, and teachers who educate student groups, just to name a portion of the potential use-cases.

First, our study shows that capturing the uncertainty and the relations between the durations of different tasks is vital to a realistic assessment of the makespan. Indeed, there are many settings in which the probabilistic knowledge about task durations is limited or costly to attain. In such circumstances, it is rather easy to design a polyhedral or ellipsoidal uncertainty set that frames the involved uncertainty. \cite{ben2009robust} provide guidance and probabilistic guarantees in favor of designing uncertainty sets that balance the level of conservatism and the probability that a constraint is violated by a scenario. Ideally, we would like to design the smallest uncertainty set that still captures the meaningful scenarios (\eg, the probability that a scenario is not included within the uncertainty set is lower than a pre-specified threshold). 

Secondly, whenever the optimal wait-and-see decisions can be taken into account in the planning stage, this should be done as it lowers the maximum possible project makespan that the scheduler can promise. In other words, a bid prepared by a decision-maker who accommodates wait-and-see decisions and thus can commit to a lower makespan (and cost) would be competitive compared to that of a bidder that does not explicitly take into account the possibility that decisions can be adapted. In particular, our experiments point out that the average advantage of adaptive-based bids is estimated to be $9-10\%$ over its non-adaptive (\ie, `regular' \gls{RO}) counterpart. \revise{We note that an adaptive policy need not necessarily be achieved by solving our \gls{MILO} formulation.  With respect to the question which adaptive policy should we use we provide two answers. The first, is obvious, use an adaptive policy which is scalable. In this context, the \gls{2SSA} heuristic is scalable and exhibits good performance, thus it could be used as an alternative to static policies. The second answer follows from our analysis of performance bounds. Specifically, we found that for large problems any rolling horizon policy that can update its the decisions at selected decision points based on the revealed information is expected to perform well -- thus, even simple time efficient policies may be enough. A good adaptive policy may be important for mid-range $\Gamma$ values rather than when uncertainty is very small (\ie,  $\Gamma\rightarrow 0$) or large (\ie, $\Gamma\rightarrow n$). Within these mid-range values, adaptive policies may provide their users a competitive edge over other policies, even the latter policies are implemented in a rolling-horizon fashion.}

While the previous point dealt with the superiority of adaptive robust policies over their static counterparts in the planning and contract stage, they are also preferable in the implementation stage. Specifically, policies that take the later-stage adaptivity of the decisions into account remain preferable even when the static policies are re-optimized every time new information becomes available (rolling-horizon). A hint into the reason for this is provided by the $42-59\%$ of the problem instances in which an adaptive policy yielded different first-stage decisions compared to a \gls{SA} policy. That means that the adaptive policies not only offer better project makespan guarantees, but also select decisions that lead to better realized duration. Our results shows clearly that the algorithms's performance deteriorate with the fraction of different firs-stage decisions compared to fully adaptive \gls{AR}. \revise{Accordingly, \gls{2SSA} (7\% different first-stage decisions) was the best performing algorithm, followed by \gls{SL} (11\%), and lastly the static \gls{SA} policy (50\%).}

A very attractive feature of the adaptive policies, as revealed through our experiments, is that their {rolling-horizon} performance is comparable to the perfect hindsight policy (\eg, the average difference between the promised and max \gls{PH} makespans was $0.0-0.1\%$ for the optimal adaptive policy compared to $5.7-9.8\%$ for the static robust policy). This suggests that the adaptive robust policy does not only protect the decision-maker against adversarial realizations of reality but it also guarantees performance close to that of the perfect hindsight policy. Thus, the typical criticism about the conservatism of static robust policies (\ie, the high price paid for robustness) does not apply to the adaptive scheduling policy.

In conclusion, while robust \gls{SA} policies are widely investigated and used in risk averse settings, they may achieve inferior performance in practice compared to adaptive alternatives. Other alternative policies such as the \gls{SL} do not scale well with the problem's size. Since the performance gap between an optimal adaptive policy and a static one is quite significant, we recommend allocating resources for finding good adaptive policies, even if those policies are not necessarily optimal. \revise{We suggested the \gls{2SSA} as such a policy that can grant its users competitive advantages both in the proposal bidding stage and in the implementation stage.}

\section{Future research} \label{sec:future.research}
Following this research, we recommend pursuing the following research directions:
\begin{enumerate}
    \item Extending the scope of scheduling problems beyond \gls{PMS}. As a first step, we suggest two problems that consider precedence constraints between  tasks: the job-shop problem, in which tasks require a unitary resource (machine) and the task includes sub-tasks with precedence constraints, and the resource constrained project scheduling problem, in which a task may require more than one type and/or more than one resource unit to be processed.
    \item It would be interesting to investigate the types of problems and settings under which first-stage decisions are different for adaptive and non-adaptive scheduling policies.
    %\item Finally, it would be worthwhile developing computationally efficient heuristic adaptive robust scheduling policies. 
\end{enumerate}

% Acknowledgments here
%\ACKNOWLEDGMENT{%
%}% Leave this (end of acknowledgment)
\section*{Acknowledgements}

We are grateful to Michael Pinedo for consultations at the early stage of this work. %Also, we thank Esther Julien for the help with the plots.

% References here (outcomment the appropriate case) 

% CASE 1: BiBTeX used to constantly update the references 
%   (while the paper is being written).
%\bibliographystyle{informs2014} % outcomment this and next line in Case 1
%\bibliography{<your bib file(s)>} % if more than one, comma separated

% CASE 2: BiBTeX used to generate mypaper.bbl (to be further fine tuned)
%\input{mypaper.bbl} % outcomment this line in Case 2
%\newpage
\bibliographystyle{apa}
\bibliography{RO_Scheduling}

\begin{APPENDIX}{}

\section{Proofs}\label{app:proofs}
\proof{Proof of Proposition~\ref{prop.box.the.same}.}
Let $J\in\mathcal{J}_{n.m}$ be some partition of the tasks to machines. Then, the makespan duration induced by this partition is given by
$$
\max_{d\in \mathcal{U}}\max_{j\in[m]}\sum_{i\in J_j} d_i= \max_{j\in[m]}\sum_{i\in J_j} \bar{d}_i,$$
and the optimal \gls{SA} policy $J^*=(J^*_1,\ldots,J^*_m)$ satisfies
$$\max_{j\in[m]}\sum_{i\in J^*_j} \bar{d}_i\leq \max_{j\in[m]}\sum_{i\in J_j} \bar{d}_i,\; \forall J\in \mathcal{J}_{n.m}\;.
$$
Now let $\bar{d}=(\bar{d}_1,\ldots,\bar{d}_n)$ and define
$J^{\text{AR}\ast}=\mathcal{M}(P^{\text{AR}\ast},\bar{d})$, where $P^{\text{AR}\ast}$ is an optimal \gls{AR} policy. Then, 
\begin{align}\label{eq:SA_equal_AR}
\max_{j\in[m]}\sum_{i\in J^*_j}\bar{d}_i\leq \max_{j\in[m]} \sum_{i\in J^{\text{AR}\ast}} \bar{d}_i\leq \max_{d\in\mathcal{U},J=\mathcal{M}(P^{\text{AR}\ast},d)} \max_{j\in[m]} \sum_{i\in J_j} d_i,
\end{align}
where the first inequality follows from the optimality of $J^*$ with respect to $J^{\text{AR}\ast}$, and the last inequality follows from adding the maximum over all $d\in \mathcal{U}$.  Since, however, the worst-case makespan of optimal \gls{AR} is always shorter than or equal to that of the optimal \gls{SA} policy, both worst-case makespans are equal, and thus both inequalities in \eqref{eq:SA_equal_AR} are in fact equalities.

We now turn to prove the equality with respect to the worst-case makespan by the optimal \gls{SL}.
We first show that for any $\pi\in \Pi_n$, using the realization $\bar{d}$ results in the worst-case makespan. Assume, to the contrary, that for some $\pi\in \Pi_n$ a worst-case makespan is achieved by a realization $\tilde{d}$ that contains a component $\pi_{j}$ such that $\tilde{d}_{\pi_{j}}<\bar{d}_{\pi_{j}}$. The starting times (and, therefore, ending times) of all tasks $\pi_i$ for $i>j$ using $\tilde{d}$ would be earlier or the same as those with $\bar{d}$, leading to a shorter or identical makespan. Thus, realization $\bar{d}$ would always lead to the worst-case makespan. Now, let $\pi^\dag\in\Pi_n$ be a permutation such that $\mathcal{M}(P^{\text{SL},\pi^\dag},\bar{d})\equiv J^*$ (such a permutation can always be constructed by the order of starting times), let $\pi^*$ be the permutation associated with the optimal \gls{SL} policy, and let $J^{\text{SL}\ast}\equiv\mathcal{M}(P^{\text{SL}\ast},\bar{d})\equiv \mathcal{M}(P^{\text{SL},\pi^*},\bar{d})$. Then, 
\begin{equation*}
\max_{d\in\mathcal{U},J=\mathcal{M}(P^{\text{SL},\pi^\dag},d)} \max_{j\in[m]} \sum_{i\in J_j} d_i
=\max_{j\in[m]}\sum_{i\in J^*_j}\bar{d}_i\leq \max_{j\in[m]} \sum_{i\in J^{\text{SL}*}} \bar{d}_i= \max_{d\in\mathcal{U},J=\mathcal{M}(P^{\text{SL},\pi^*},{d})} \max_{j\in[m]} \sum_{i\in I_j} d_i,
\end{equation*}
where the first equality follows from the definition of $\pi^\dag$, \revise{the inequality follows from the optimality of the static policy $J^*$ (from the point of view of \gls{SA}, $J^{\text{SL}\ast}$ is suboptimal), the second equality follows from the optimality of $\bar{d}$ for the worst-case,} and the last equality from the definition of $J^{\text{SL}*}$ and again the fact that $\bar{d}$ is a worst case for $\pi^*$.
Combining this inequality with the optimality of $\pi^*$, we obtain that the worst-case makespans of the optimal \gls{SL} and \gls{SA} policies must be equal.\Halmos
\endproof

\revise{\proof{Proof of Theorem~\ref{theo:SA}.}
In order to lower bound the worst case makespan of the PH policy, we can restrict the adversary to split the uncertainty equally between the tasks. For these scenarios the PH policy would use allocation $\tilde{x}$. Thus,
\begin{align*}
    PH&=\max_{d\in U}\min_{x\in\{0,1\}^n} \max \{ x^\top d, (e - x)^\top d \}\\
    &\geq \max \{ \tilde{x}^\top  (d^0+\frac{\Gamma}{n}\bar{d}), (e - \tilde{x})^\top  (d^0+\frac{\Gamma}{n}\bar{d}) \}\\
    &=\max\left\{\sum_{i\in \mathcal{I}} d^0_i(1+\frac{\alpha\Gamma}{n}),\sum_{i\in \bar{\mathcal{I}}} d^0_i(1+\frac{\alpha\Gamma}{n})\right\}\\
    &=(1+\frac{\alpha\Gamma}{n})\max\left\{\sum_{i\in\mathcal{I}} d^0_i,\sum_{i\in\bar{\mathcal{I}}} d^0_i\right\}
\end{align*}
Moreover, the SA worst case is not worse than the one obtain by using $\tilde{x}$ as the allocation, and thus,
\begin{align*}
    SA&= \min_{x\in\{0,1\}^n}\max_{d\in U} \max \{ x^\top d, (e - x)^\top d \}\\
    &\leq  \max_{d\in U}\max \{ \tilde{x}^\top d, (e - \tilde{x})^\top d \}\\
    &= \max_{d\in U} \max\left\{\sum_{i\in\mathcal{I}} d_i,\sum_{i\in \bar{\mathcal{I}}} d_i\right\}\\
    &=\max\left\{\sum_{i\in\mathcal{I}} d^0_i+\alpha\max_{\substack{\mathcal{S}\subseteq\mathcal{I}, |\mathcal{S}|=\min\{\lfloor\Gamma\rfloor,|\mathcal{I}|\}\\
    j\in \mathcal{I}\setminus\mathcal{S}}}\left\{ \sum_{k\in\mathcal{S}} d^0_k+(\min\{\Gamma,|\mathcal{I}|\}-|\mathcal{S}|)d^0_j\right\},\right.\\
    &\; \left.\sum_{i\in [n]\setminus \mathcal{I}} d_i^0+\alpha\max_{\substack{\mathcal{S}\subseteq[n]\setminus \mathcal{I}, |\mathcal{S}|=\min\{\lfloor\Gamma\rfloor,n-|\mathcal{I}|\}\\
    j\in [n]\setminus \mathcal{I}\cup\mathcal{S}}} \left\{\sum_{k\in\mathcal{S}} d^0_k+(\min\{\Gamma,n-|\mathcal{I}|\}-|\mathcal{S}|)d^0_j\right\}\right\}.
    %&\leq \max\{\sum_{i\in\mathcal{I}} d^0_i,\sum_{i\in[n]\setminus \mathcal{I}} d^0_i\}+\max_{\mathcal{S}\subset[n],|\mathcal{S}|=\lfloor\Gamma\rfloor, j\in[n]\setminus \mathcal{S}} \sum_{k\in\mathcal{S}}\alpha d^0_k+(\Gamma-\lfloor\Gamma\rfloor)\alpha d^0_j.
\end{align*}
Using our defined notation, can rewrite the above inequality as
\begin{align*}
SA&\leq \max\left\{\sum_{i\in\mathcal{I}} d^0_i,\sum_{i\in\bar{\mathcal{I}}} d^0_i\right\}+\alpha\max_{\mathcal{S}\in \{\mathcal{I},\bar{\mathcal{I}}\}}\left\{\sum_{k=1}^{\min\{|\mathcal{S}|,\lfloor{\Gamma}\rfloor\}} d^{0,\mathcal{S}}_{(k)}+\Delta^\mathcal{S}(\Gamma)d^{0,\mathcal{S}}_{\left(\min\{|\mathcal{S}|,\lfloor\Gamma\rfloor\}+1\right)}\right\}. \end{align*}
Dividing the two bounds we obtain
\begin{align}\label{eq:ub}
\frac{\text{SA}}{\text{PH}}&\leq \frac{1}{1+\frac{\alpha\Gamma}{n}}+\frac{\alpha\max_{\mathcal{S}\in \{\mathcal{I},\bar{\mathcal{I}}\}}\left\{\sum_{k=1}^{\min\{|\mathcal{S}|,\lfloor{\Gamma}\rfloor\}} d^{0,\mathcal{S}}_{(k)}+\Delta^\mathcal{S}(\Gamma)d^{0,\mathcal{S}}_{\left(\min\{|\mathcal{S}|,\lfloor\Gamma\rfloor\}+1\right)}\right\}}{(1+\frac{\alpha\Gamma}{n})\max\left\{\sum_{i\in \mathcal{I}} d^0_i,\sum_{i\in [n]\setminus \mathcal{I}} d^0_i\right\}}\nonumber\\
&\leq \frac{n}{n+\Gamma\alpha}+\frac{\alpha n}{n+\alpha\Gamma}\max_{\mathcal{S}\in \{\mathcal{I},\bar{\mathcal{I}}\}}\frac{\sum_{k=1}^{\min\{|\mathcal{S}|,\lfloor{\Gamma}\rfloor\}} d^{0,\mathcal{S}}_{(k)}+\Delta^\mathcal{S}(\Gamma)d^{0,\mathcal{S}}_{\left(\min\{|\mathcal{S}|,\lfloor\Gamma\rfloor\}+1\right)}}{\sum_{i\in\mathcal{S}} d_i^0} \nonumber \\ 
&=\frac{n}{n+\Gamma\alpha} \left( 1+\alpha\max_{\mathcal{S}\in \{\mathcal{I},\bar{\mathcal{I}}\}}\frac{\sum_{k=1}^{\min\{|\mathcal{S}|,\lfloor{\Gamma}\rfloor\}} d^{0,\mathcal{S}}_{(k)}+\Delta^\mathcal{S}(\Gamma)d^{0,\mathcal{S}}_{\left(\min\{|\mathcal{S}|,\lfloor\Gamma\rfloor\}+1\right)}}{\sum_{i\in\mathcal{S}} d_i^0} \right).
\end{align}}

\revise{We now look at  \eqref{eq:ub} for different values of $\Gamma$. Specifically, we look at the term inside the maximum. When $\Gamma=0$, $\Delta^\mathcal{S}(\Gamma)=0$, and thus it is equals zero and upper bound is equal to $1$. When $\Gamma=n$ then again $\Delta^\mathcal{S}(\Gamma)=0$, its numerator and denominator are equal, and so the upper bound is again equal to $1$. 
\endproof }

\revise{\proof{Proof of Theorem~\ref{theo:RH}.}
By construction we have that
\begin{equation}\text{PH}\geq \max_{d\in U} \frac{\sum_{i=1}^n d_i}{2}.\end{equation}
Moreover, as explained above, for any RH methodology, we have the following upper bound
\begin{equation}\text{RH}\leq \max_{d\in U, j\in[n]} \frac{\sum_{i=1}^n d_i}{2}+\frac{d_j}{2}.\end{equation}
Using these bounds and Assumption~\ref{ass:shape_uncertainty}, we have that
\begin{align}
 \frac{\text{RH}}{\text{PH}}&\leq\frac{\max\limits_{d\in U,j\in[n]}\left\{\sum\limits_{i=1}^n \frac{d_i}{2}+ \frac{d_j}{2}\right\}}{\max\limits_{d\in U}\sum\limits_{i=1}^n \frac{d_i}{2}}\\
 &\leq \frac{\max_{d\in U}\sum\limits_{i=1}^{n} \frac{d_i}{2}+\frac{\bar{a}\min\{(1+\alpha),(1+\alpha\Gamma)\}}{2}}{\max_{d\in U}\sum\limits_{i=1}^{n} \frac{d_i}{2}}\\
 &\leq 1+\frac{\bar{a}\min\{(1+\alpha),(1+\alpha\Gamma)\}}{\ubar{a}(n+\alpha\Gamma)}.
\end{align}
\endproof}

\section{Extensions for more than two machines}
\subsection{The scheduler's DP}\label{app:DP_more_than_two}
When  $m>2$, more than one task may still be in process when a given task completes. Thus, we extend the definition of a state to 
$S,F,D,I,\bar{D}$, where $I$ is the set of indices of running tasks, and their respective processing duration thus far is represented by vector $\bar{D}$. Decisions can be made when there is an idle machine, \ie, when one task completes ($|I|<m$,) and there are still tasks to schedule, ($|S|<n$). When $|S|=n$, all decisions have been already made.

Similar to $m=2$ machines, at each decision point, the scheduler seeks the allocation policy that achieves the minimal worst-case remaining makespan. The recursive formulation is:
%{\small
\begin{align}\label{eq:T_no_Waitng_eps_m_more_than_2}
&T(S, F,D,I,\bar{D})=\\
&\quad \min_{k\notin S}%F\cup I}
\max\Biggl\{ \max_{\substack{d_k:d\in U_{S,F,D,I, \bar D},\\ d_k \leq \min\limits_{j \in I} (d_j-\bar{D}_{j})}} d_k+
T([S,k], [F,k],[D,d_k],I, (\bar{D}_{j}+d_k)_{j\in I} )\nonumber\\
&\qquad \max_{\substack{(l,d_l):\\d\in U_{S,F,D,I, \bar D}\\l=\argmin\limits_{j \in I} (d_j-\bar{D}_{j})\\d_k>  d_l-\bar{D}_{l} }}  d_l-\bar{D}_{l}+
T([S,k],[F{,}l],[D{,}d_l],[I\setminus \{l \}, k],[(\bar{D}_{j})_{j {\in} I\setminus\{l\}},0]+d_l{-}\bar{D}_{l}) \Biggr\}.\nonumber
\end{align} %}
The objective is to allocate the task $k$ that minimizes the worst-case remaining makespan. For each $k$, there are two types of worst-case scenarios to consider. The first term in the first $\max$ of  \eqref{eq:T_no_Waitng_eps_m_more_than_2} relates to the possibility that under the worst-case scenario, $k$ completes its processing before the completion of any of the other tasks already underway in $I$. In such a case, the next decision point is when task $k$ completes and all tasks in $I$ are still in processing (thus, the composition of $I$ does not change). The second term describes the scenarios in which the next task completes before the newly scheduled task $k$. Thus, the time until the next decision point is $\min_{j \in I} (d_j-D_{j})$, in which case the composition of set $I$ changes to include $k$ and to exclude the task that just finished processing.

The boundary case occurs when all tasks have been scheduled, yet some of them are still in processing; that is, $|S|=n$, $|F|<n$. In such a case, the remaining makespan amounts to the time until the completion of the last task among the tasks still being processed:
\begin{align}
&T(S, F,D,I,\bar{D})=\max_{i\in I} \max_{d_i:d\in U_{S, F,D,I,\bar{D}}}d_i-\bar{D}_i. \label{eq:T_no_Waitng_eps_m_more_than_2_boundary} 
\end{align}

\subsection{The adversary's DP}\label{app:DP_more_than_two_advers}
As in the two-machine formulation, the adversary makes her decisions immediately after a new task $k$ has been scheduled in states $([S,k],F,D,I,\bar{D})$ where $k\notin S\equiv F\cup I$. The recursive formula for the worst-case remaining makespan $T'$ is:
{\small
\begin{align}\label{eq:adversary_DP_m_machines}
\begin{split}
&T'([S,k], F,D,I,\bar{D})=\max\Bigl\{\Bigr. \max_{\substack{d_k:d\in U_{[S,k],F,D,I,\bar D},\\d_k \leq \min_{j \in I} (d_j-\bar{D}_{j})}} d_k+\min_{l \notin S\cup\{k\}} T'([S,k,l], [F,k],[D,d_k],I,\bar D+d_k\}, 0]),\\
&\quad \max_{\substack{(i,d_i):d\in U_{[S,k],F,D,I,\bar D},\\i\in\argmin_{j \in I} (d_j-\bar{D}_{j}), \\ d_k>  d_i-\bar{D}_i}} d_i-\bar{D}_{i}+ \min_{l \notin S\cup\{k\}} T'([S,k,l],[F{,}i],[D{,}d_i],[I\setminus\{i\},k],[(\bar{D}_{j}+d_i{-}\bar{D}_{i})_{j\in I\setminus\{i\}}, d_i{-}\bar{D}_{i}] \Bigl.\Bigr\}, 
\end{split}
\end{align}}
as long as $|S|\leqslant n-2$ (meaning that there are still more tasks to schedule). For the case $|S|=n-1$ and $|F|<n$, we have
\begin{align*}
& T'([S,k], F,D,I,\bar{D}]) = \max\left\{\max_{i\in I} \max_{d_i:d\in U_{[S,k], F,D,I,\bar{D}}} d_i-\bar{D}_i,\max_{d_k:d\in U_{[S,k], F,D,I,\bar{D}}} d_k\right\}
\end{align*}
which is the longest possible time until the last task completes.

\subsection{MIO formulation}\label{More_than_two_MIO}

We consider a setting with $m$ machines and $n$ tasks ($m<n$) and denote the set of all states $\sigma=(S,F)$ as $\mathcal{V}$. To decrease the state space of the scenario tree, certain states are discarded (see also the example described in Section~\ref{sec:2machines4tasks}). In particular,
\begin{itemize}
\item we omit states where the length of $S$ is smaller than $m$, since they correspond to an initialization of the system in which at least one machine is idle, thus additional tasks must be scheduled immediately. Hence, the states that follow the initial state are characterized by $|S|=m$ (all machines are busy).
\item Similarly, the last decision to be made is the one after which a single task remains to be scheduled. Such a decision is made at a state where the length of $S$ is $n - 2$ and the length of $F$ is $n - m - 1$. Afterwards, the adversary decides which tasks to complete and in which order. Hence, the states that follow states with  $|S|=n-1$ and $|F|=n - m - 1$ are the end states in which $|S|=|F|=n$.
\end{itemize}
Given this, we consider the state space $\mathcal{V}$ consisting of
the sets of the scheduler states $\mathcal{D}$ and adversary states $\mathcal{N}$ (Table \ref{tab:scenario_states_schedular}).
\begin{center}
\begin{table}[ht]
\caption{The sets of states for the scheduler $\mathcal{D}$ (left) and the adversary $\mathcal{N}$ (right).} \label{tab:scenario_states_schedular}
\vspace{10pt}
\begin{minipage}{.45\textwidth}
{\centering
{\def\arraystretch{1.2}\tabcolsep=3pt
{\scriptsize 
\begin{tabular}{p {0.45\textwidth} p {0.45\textwidth}} 
\hline
\textbf{$\boldsymbol{S}$-length} & \textbf{$\boldsymbol{F}$-length} \\ \hline
0 & 0 \\ \hline
$m$ & 1 \\\hline
$m + 1$ & 2 \\\hline
\vdots & \vdots \\\hline
$n - 2$ & $n - m - 1$\\ \hline
\end{tabular}}}}\\
\end{minipage}
\begin{minipage}{.45\textwidth}
{\centering
{\def\arraystretch{1.2}\tabcolsep=3pt
{\scriptsize 
\begin{tabular}{p {0.45\textwidth} p {0.45\textwidth}} 
\hline
\textbf{$\boldsymbol{S}$-length} & \textbf{$\boldsymbol{F}$-length} \\ \hline
$m$ & 0 \\ \hline
$m + 1$ & 1 \\\hline
\vdots & \vdots \\\hline
$n - 1$ & $n - m - 1$\\ \hline
$n$ & $n$  \\
\hline
\end{tabular}}}}\\
\end{minipage}
\end{table}
\end{center}
With these formulations, we can define the \gls{MILO} formulation in the same way as in Section~\ref{sec:MIO}.

\section{Complexity of the \gls{MILO} formulation (\ref{opt-P})} \label{appendix:milo.complexity}
In this appendix, we study the complexity of the \gls{MILO} formulation \eqref{opt-P} needed to solve the scheduling problem. The scheduler's first decision concerns which $m$ tasks out of the $n$ to schedule, followed by the adversary's decision regarding which of the $m$ tasks to end first. Next, the scheduler has to choose one of the remaining tasks to schedule to the newly idle machine. This sequence repeats itself until the scheduler has no more tasks to schedule, in which case the adversary has to decide on the termination order of the $m$ in-process tasks. Thus, the number of end states in the scenario tree of the mixed-integer formulation is given by
\begin{align*}
|\mathcal{L}| = & {\binom{n}{m}}  m  (n - m )  m  (n - m - 1)\cdot  \ldots  \cdot m \cdot 2 \cdot m! \\
= & {\binom{n}{m}} {m^{n - m - 1}} (n - m )!  m! \\
= & {n!} m^{n - m - 1}.
\end{align*}
The number of scheduler nodes is:
\begin{align*}
|\mathcal{D}| = & 1  + {\binom{n}{m}}  m  + {\binom{n}{m}}  m  (n - m )  m   + \ldots 
 + {\binom{n}{m}}  m  (n - m )  m  (n - m -1)  m  \ldots 3\cdot m \\
& =1+\sum_{i=0}^{n-m-2} {\binom{n}{m}}{\frac{(n-m)!}{(n-m-i)!}}m^{i+1}\\
& =1+\frac{n!}{m!} \sum_{i=1}^{n-m-1}\frac{m^{i}}{(n-m-i+1)!}.
\end{align*}
The number of binary variables required for the nature nodes' maximum reformulation as \gls{MILO} constraints is, therefore, $|\mathcal{D}|+|\mathcal{L}|-1$. Thus, the \gls{MILO} formulation requires an exponential number of binary variables.

Consequently, the \gls{MILO} problem \eqref{opt-P} is not scalable. In the next section we conduct a numerical study with small problems to compare the \gls{MILO} formulation (\ie, \gls{AR}) to other alternatives (\eg,  \gls{SA}, \gls{SL}), which will enable us to evaluate the possible benefits of applying an adaptive or at least partially adaptive policy for which finding solutions is less computationally demanding.

\section{\revise{Column-and-constraint generation algorithm for the 2SSA}} \label{app:ccg}
\revise{Assume we schedule two tasks $i$, $j$ as the first ones, and our goal is to compute the worst-case makespan of the heuristic for the case $d_i < d_j$ (first term in the outer $\max$ in \eqref{eq:two.stage.SA}). Next, proceed as follows
\begin{enumerate}
    \item Initialize a list of static allocation schedules $\mathcal{X} = \{x^1\}$ consisting of an arbitrary $x^1 \in \{ 0, 1 \}^n$, where $x^1_i = 1$, and $x^1_j = 0$.
    \item Solve the master problem
    \begin{align}
    \max\limits_{t, t^{k}, d^{k} \in U} \ & t \label{eq:ccg.masterproblem} \\
    \text{s.t.} \ & t \leq t^{k} && \forall k = 1\ldots, |\mathcal{X}| \nonumber \nonumber \\
    & t^k \leq \max \{ x^{k\top} d^{k}, (\mathbf{e} - x^{k})^\top d^{k} \} && \forall k = 1\ldots, |\mathcal{X}| \nonumber \\
    & d^k_i = d^l_i && \forall k = 1\ldots, |\mathcal{X}| \nonumber \\
    & d^k_i \leq d^k_j. && \forall k,\nonumber 
    \end{align}
    Denote by $\bar{d}_i$ the optimal value of $d_i$, and by $\bar{t}$ the optimal worst-case makespan. Go to Step 3,
    \item Solve the subproblem:
    \begin{align}
        \min\limits_{x\in \{0, 1 \}^n, v} \ & v \label{eq:ccg.subproblem} \\ 
        \text{s.t.} \ & v \geq \max\left\{ \sup\limits_{d \in U: \ \bar{d}_i = d_i \leq d_j} x^\top d, \sup\limits_{d \in U: \ \bar{d}_i = d_i \leq d_j} (\mathbf{e} - x)^\top d \right\}\label{eq:subproblem.objective} \\
        & x_i = 1 \nonumber \\
        & x_j = 0 \nonumber
    \end{align}
If the optimal value $\bar{v} < \bar{t}$, set: $\mathcal{X} := \mathcal{X} \cup \bar{x}$ and go to Step 2. Otherwise, finish. 
\end{enumerate}}

\revise{The interpretation of the master problem \eqref{eq:ccg.masterproblem} is as follows. The adversary is seeking a scenario tree of task durations to maximize the worst-case makespan across all second-stage schduler's decision possibilities $x^k$. Because for each selected $x^k$ the worst-case realization of $d$ might be different, for each $x^k$ we have a separate decision vector $d^k$. However, these can differ in all tasks except for task $i$ because of the non-anticipativity of the adversary. The first constraint means that the scheduler will pick the schedule that minimizes the worst-case duration, from the available options.}

\revise{Subproblem \ref{eq:ccg.subproblem} answers the following question: for the event $\bar{d}_i = d_i < d_j$, does there exist a static allocation schedule $x \in \{0, 1 \}^n$, $x_i = 1$, $x_j = 0$, not included in $\mathcal{X}$, and which would make the worst-case makespan better? If the answer is `yes', the decision possibility is added and the master problem is re-solved, until convergence. Importantly, the maximum term issue in \eqref{eq:subproblem.objective} can be reformulated using the big-M trick, and if the set $U$ is a polytope, the inner supremizations in \eqref{eq:subproblem.objective} can be dualized so that the final problem is an MILP.}

\section{Box uncertainty set sampled results}\label{app:box_uncertainty_sampled}

\begin{table}[ht!]
\TABLE
{Performance measure results for type II uncertainty sets. \label{tab:results.box}}
{\centering
{\def\arraystretch{1.2}\tabcolsep=3pt
{\scriptsize \begin{tabular}{lccccc} \hline
\textbf{Short name} & \textbf{\gls{AR}} & \textbf{\gls{2SSA}} & \textbf{\gls{SL}} & \textbf{\gls{SA}} & \textbf{\gls{PH}}\\  \hline
\revise{Suboptimal initial} & - & \revise{$2\%$} & \revise{$4\%$} & \revise{$32\%$} & - \\ \hline
Promised makespan & $ 7.784 (1.430) $ & $ 7.788 (1.428) $ & $ 7.798 (1.429) $ & $ 8.217 (1.470) $ & - \\ \hline 
Max makespan & $ 7.784 (1.430) $ & $ 7.787 (1.429) $  & $ 7.791 (1.429) $ & $ 7.913 (1.425) $ & $7.783 (1.430)$ \\ \hline 
%Max makespan non RH & $ 7.78 $ & $ 7.787 $  & $ 7.798 $ & $ 8.21 $ & - \\ \hline 
Makespan & $ 6.105 (1.493) $ & $ 6.097 (1.496) $  & $ 6.098 (1.495) $ & $ 6.139 (1.513) $ &$ 5.966 (1.493)$ \\ \hline 
%Makespan non-RH & $ 6.105 $ & $ 6.21 $  & $ 6.22 $ & $ 6.60 $ & - \\ \hline 
\end{tabular}}}}
{Makespan values are presented with their standard deviations in parentheses.} %Ratio values represent the difference, in percentage, with respect to the value of the policy in the denominator of the ratio and their standard deviations in parentheses.}
\end{table}

\begin{figure}
\centering
\includegraphics[width=0.9\textwidth]{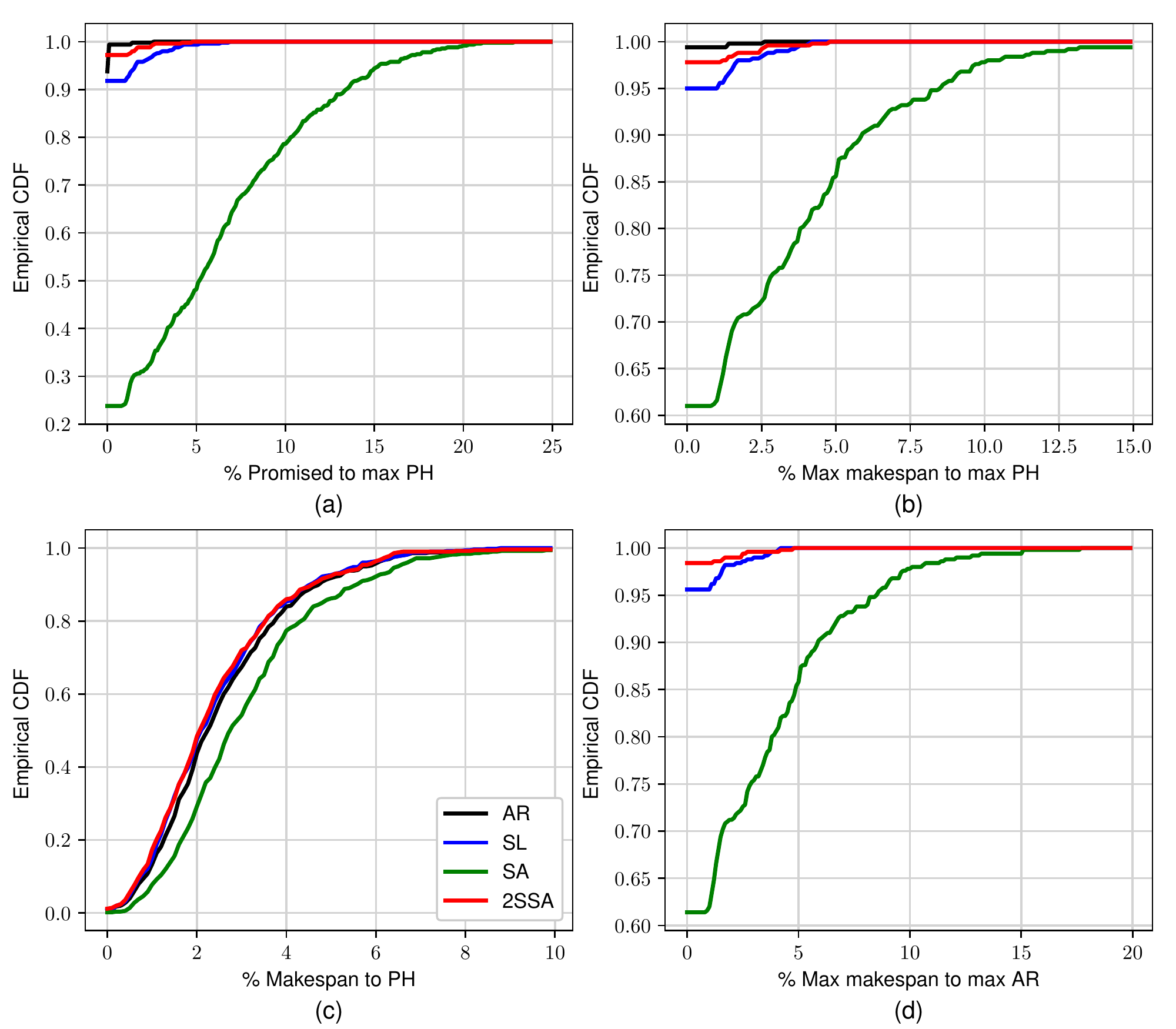}
\caption{\revise{Empirical cumulative distribution functions (over instances) of the last four performance measures of Table~\ref{tab:results.measures} for Type II uncertainty set.}}
\label{fig:ecdfs_type_0}
\end{figure}

For \gls{SA}, \gls{SL}, and \gls{2SSA}, the percentages of different first-stage decisions, compared to \gls{AR}, are  $32\%$, $4\%$, and $2\%$ respectively.

The \gls{AR} promised makespan was the shortest, with \gls{2SSA} being nearly the same, \gls{SL} very close behind (longer by $0.17\%$), and \gls{SA} trailing behind with a longer makespan by $5.3\%$. The upper bound on the ``price'' that a risk-averse scheduler, who commits to the promised makespan, is expected to pay upfront with respect to a \gls{PH} policy (as indicated by the Promised to max \gls{PH} measure) was very small for \gls{AR} ($0.0\%$), \gls{2SSA} ($0.1\%$) \gls{SL} ($0.2\%$), and much higher for \gls{SA} ($5.8\%$).

\end{APPENDIX}
%\newpage
%\printglossary
\end{document}